\documentclass[12pt]{amsart}
\usepackage{epsfig,amsthm,amsmath,amsfonts,amssymb,bbm,
color,nicefrac,comment,hyperref}
\usepackage[normalem]{ulem}
\usepackage{tikz}
\usepackage[export]{adjustbox}
\usepackage{tikz,pgfplots,pgfplotstable}
\usetikzlibrary{arrows,decorations,backgrounds,positioning}

\usepackage[bb=boondox]{mathalfa} 
\newtheorem{theorem}{Theorem}

\numberwithin{equation}{section}

\theoremstyle{definition}
\newtheorem{definition}[theorem]{Definition}
\theoremstyle{remark}
\newtheorem{remark}[theorem]{Remark}

\newcommand{\E}{\mathbb{E}}
\newcommand{\V}{\mathbb{V}}
\newcommand{\Cor}{\mathbb{C}\!\operatorname{or}}
\newcommand{\Cov}{\mathbb{C}\!\operatorname{ov}}

\renewcommand{\d}{\operatorname{d}\!}
\newcommand{\Rset}{\mathbb{R}}

\newcommand{\bfXi}{\boldsymbol{\Xi}}
\newcommand{\bfpsi}{\boldsymbol{\psi}}
\newcommand{\bfchi}{\boldsymbol{\chi}}
\newcommand{\bfp}{{\bs p}}
\newcommand{\bfx}{{\bs x}}
\newcommand{\bfC}{{\bs C}}
\newcommand{\bfL}{{\bs L}}
\newcommand{\bfM}{{\bs M}}
\newcommand{\bfT}{{\bs T}}
\renewcommand{\S}{\mathbb{S}}
\newcommand{\bbS}{\mathbbb{S}}

\newcommand{\isdef}{\mathrel{\mathrel{\mathop:}=}}

\newcommand{\refd}{{\operatorname{ref}}}
\newcommand{\bs}[1]{{\boldsymbol#1}}
\newcommand{\xref}{\widehat{\bs x}}
\newcommand{\Dtrace}{\gamma_0^{\mathrm{int}}}
\pgfplotsset{compat=1.16}
\begin{document}

\title[Multilevel quadrature in acoustic scattering]
{Isogeometric multilevel quadrature for forward and inverse random acoustic scattering}
\author[J.~D\"olz]{J\"urgen D\"olz}
\address{J\"urgen D\"olz,
Institute for Numerical Simulation,
University of Bonn, 
Endenicher Allee 19b, 53115 Bonn, Germany.}
\email{doelz@ins.uni-bonn.de}
\author[H.~Harbrecht]{Helmut Harbrecht}
\address{Helmut Harbrecht, 
Department of Mathematics and Computer Science,
University of Basel, 
Spiegelgasse 1, 4051 Basel, Schweiz.}
\email{helmut.harbrecht@unibas.ch}
\author[C. Jerez-Hanckes]{Carlos Jerez-Hanckes}
\address{
Carlos Jerez-Hanckes,
Faculty of Engineering and Sciences, Universidad Adolfo Ib\'a\~{n}ez, Santiago, Chile.}
\email{carlos.jerez@uai.cl}
\author[M.~Multerer]{Michael Multerer}
\address{
Michael Multerer,
Institute of Computational Science,
USI Lugano,
Via Giuseppe Buffi 13, 6900 Lugano, Schweiz}
\email{michael.multerer@usi.ch}
\date{}
\begin{abstract}
We study the numerical solution of forward and inverse 
acoustic scattering problems by randomly shaped obstacles 
in three-dimensional space using a fast isogeometric boundary
element method. Within the isogeometric
framework, realizations of the random scatterer 
can efficiently be computed by simply updating the NURBS
mappings which represent the scatterer. This way, we
end up with a random deformation field.
In particular, we show that the knowledge of the
deformation field's expectation and covariance at the surface
of the scatterer are already sufficient to compute
the surface Karhunen-Lo\`eve expansion. Leveraging on the
isogeometric framework, we utilize
multilevel quadrature methods for the efficient
approximation of quantities of interest, such as the
scattered wave's expectation and variance.
Computing the wave's Cauchy data at
an artificial, fixed interface enclosing the random obstacle,
we can also directly infer quantities of interest in free space.
Adopting the Bayesian paradigm, we finally compute the expected
shape and the variance of the scatterer from noisy measurements
of the scattered wave at the artificial interface.
Numerical results for the forward and
inverse problem are given to demonstrate the feasibility 
of the proposed approach.
\end{abstract}
\maketitle 

\section{Introduction}\label{sec:intro}
The reliable computer simulation of phenomena where
acoustic waves are scattered by obstacles is of great importance 
in many applications. These include for example the modelling of 
sonar and other methods of acoustic location, as well as outdoor 
noise propagation and control, especially stemming from automobiles, 
railways or aircrafts. Since an analytical solution of scattering problems 
is in general impossible, numerical approaches are called for. 

Most acoustic scattering problems may be formulated in the 
frequency domain by employing the Helmholtz equation: assume an 
acoustic wave encounters an impenetrable, bounded obstacle $D\subset
\mathbb{R}^3$, having a Lipschitz smooth boundary $S
\isdef\partial D$, and, as a consequence, gets scattered. Then, describing the incident plane wave $u_{\text{inc}}({\bs x}) = e^{i\kappa\langle{\bs d},
{\bs x}\rangle}$ with known wavenumber $\kappa$ and direction ${\bs d}$, where $\|{\bs d}\|_2=1$, the total wave 
\[
  u = u_{\text{inc}}+u_{\text{s}}
\]
is obtained by solving the exterior boundary value problem
\begin{equation}\label{eq:pde}
 \begin{aligned}
  \Delta u + \kappa^2 u = 0\quad&\text{in}\ \mathbb{R}^3\setminus\overline{D},\\
  	u = 0\quad&\text{on}\ S,\\
	  \sqrt{r}\bigg(\frac{\partial u_{\mathrm{s}}}{\partial r}-i\kappa u_{\mathrm{s}}\bigg)
  	\to 0\quad&\text{as}\ r = \|{\bs x}\|_2\to\infty.
  \end{aligned}
\end{equation}
The homogeneous Dirichlet condition at $S$ corresponds 
to a \emph{sound-soft} obstacle, whereas a homogeneous Neumann 
condition would correspond to a \emph{sound-hard} obstacle.
The function $u_{\mathrm{s}} = u-u_{\mathrm{inc}}$ is called the \emph{scattered wave}.
Although we restrict ourselves here to the
sound-soft case, the presented concepts are also suitable 
to treat sound-hard obstacles as well as for penetrable obstacles, i.e.~objects described by a different 
diffractive index to the free space.

In this article, we consider the situation of a randomly shaped obstacle $D=D({\bs y})$,
where \({\bs y}\in\Gamma\subset\mathbb{R}^\mathbb{N}\) is some random parameter. 
This shape uncertainty might for example issue from measurement or modelling errors. 
As a consequence, the total wave itself becomes a random field $u({\bs y})$. 
Our goal is to compute  the first and second order statistics of the scattered 
wave, these are the expectation $\E[u_{\mathrm{s}}]$ and the variance
$\V[u_{\mathrm{s}}]$. We especially demonstrate how to compute 
the scattered wave's second moment in a deterministic 
fashion from its Cauchy data's second moment on an 
artificial, fixed interface $T$, which almost surely 
encloses the domain $D({\bs y})$. In combination with low-rank 
techniques, this drastically reduces the high dimensionality 
of the random scattering problem, compare \cite{HIM}.
In order to speed up the computations of the Cauchy data's
statistics even further, we employ the multilevel quadrature method,
see e.g.\ \cite{BSZ11,Gil15,H2,HPS12}.

Our approach lies in the {\em domain mapping} category as it transfers the shape uncertainty onto a fixed reference domain, and allows to deal with large deformations (see \cite{JSZ17,AJZ20}). In contrast, {\em perturbation techniques} resort to shape derivatives to linearize fields for small deviations with respect to both wavelength and scatterers' shape from a nominal, reference geometry. By Hadamard's theorem, the resulting linearized equations for the first order shape sensitivities are homogeneous equations posed on the nominal geometry, with inhomogeneous boundary data only . Using first-order shape Taylor expansions,
one can derive tensor deterministic first-kind boundary integral equations for the statistical moments of
the scattering problems considered. These are then approximated by sparse tensor Galerkin discretizations
via the combination technique (\emph{cf.}~\cite{EJH20} and references therein). Though successfully applied to three-dimensional Helmholtz Dirichlet, Neumann, impedance and transmission problems \cite{EJH20}, and even for diffraction gratings \cite{SAF18}, therein random perturbations are required to be sufficiently small. High-order approaches \cite{HDE18,Dol2020} lead to at least third order accurate approximations with respect to the perturbation amplitude of the domain variations. Finally, in \cite{castrillonc2017hybrid} a hybrid between domain-mapping and perturbation methods was presented.

For the numerical realization of random obstacles, we employ the 
methodology from \emph{isogeometric analysis} (IGA). IGA has been
introduced in \cite{HCB05} in order to incorporate simulation techniques 
into the design workflow of industrial development and thus allows to deal
with domain deformations in a straightforward manner. By representing the
geometry and domain deformations by \emph{non-uniform rational B-splines}
(NURBS), realizations of the random scatterer can efficiently be computed
by simply updating the NURBS mappings which represent the scatterer.
In addition, the naturally emerging sequence of nested approximation spaces
can directly be employed in multilevel quadrature methods. With
regard to the isogeometric boundary element approach for the
scattered wave computations, compare \cite{DHK+18,FGHP17,SBTR12,TM12}, 
we show that all computations can directly be performed at the boundary of
the deformed scatterer. This particularly applies to
the random deformation field which only needs to
be computed with respect to a reference surface.
This way, we can model large deformations without having to deal with
very fine volume meshes, which would otherwise be necessary to
properly resolve the deformation field within the scatterer.
Moreover, the meshing of the unbounded free space is avoided. 
Therefore, the isogeometric boundary element method is the method
of choice for the problem at hand.
For the numerical computations, we rely on
the fast isogeometric boundary element method developed in \cite{DHK+20,DHK+18,DHP16,
DKSW,HP13}, which is available as \verb|C++|-library 
\verb+bembel+ \cite{bembel, DHK+20}. In order to speed up computations,
\verb+bembel+ utilizes $\mathcal{H}^2$-matrices with the interpolation based
fast multipole method \cite{GR87,GR91,HB02}.

To our knowledge, the present work constitutes the first fast IGA implementation for time-harmonic acoustic wave scattering for shape uncertainty quantification. Having this fast forward solver at our disposal, we
also consider acoustic shape inversion by Bayesian inference: 
Given noisy measurements of the scattered wave at certain locations in
free space, we determine statistics of the uncertain scatterer's
shape. To this end, we employ the multilevel ratio estimator,
see \cite{DGLS17} and the references therein, and compute the expected shape and its variance.

The rest of article is organized as follows: In Section~\ref{sec:randomDomains} is concerned with the modeling of random domains and their parametrization by means of a Karhunen-Lo\`eve expansion. Section~\ref{sec:DiscRD} we perform the efficient discretization of the random deformation field by means of isogeometric analysis. In Section~\ref{sec:bie} we introduce the boundary integral formulation of the problem under consideration and discuss the use of the artifical interface for the representation of the scattered wave and its statistics. Section~\ref{sec:MLQ} briefly recalls the multilevel quadrature method, whereas Section~\ref{sec:bayes} recalls its application to Bayesian inference. Finally, Section~\ref{sec:numex} is devoted to numerical examples showcasing the ideas discussed.

\section{Random domain model}\label{sec:randomDomains}
\subsection{Modelling of random domains}
In what follows, let $D_\refd\subset\mathbb{R}^3$ denote a 
Lipschitz domain with piecewise smooth surface \(S_\refd\isdef
\partial D_\refd\) and let \((\Omega,\mathcal{F},\mathbb{P})\) be 
a complete probability space. We assume that the uncertainty 
in the obstacle is encoded by a random deformation field, cf.\ 
\cite{HPS16}. We hence assume the existence of a uniform 
\(C^1\)-diffeomorphism \({\bs\chi}_{D}\colon\overline{D_\refd}
\times\Omega\to\mathbb{R}^3\), i.e.\
\begin{equation}\label{eq:unifVfield}
\|{\bs\chi_{D}}(\omega)\|_{C^1(\overline{D_\refd};\mathbb{R}^3)},
\|{\bs\chi}_{D}^{-1}(\omega)\|_{C^1(\overline{D(\omega)};\mathbb{R}^3)}
\leq C_{\operatorname{uni}}\quad\text{for\(\mathbb{P}\)-a.e.\ }\omega\in\Omega,
\end{equation}
such that 
\[
D(\omega)={\bs\chi}_{D}(D_\refd,\omega).
\]
Particularly, since \({\bs\chi}_{D}\in L^\infty\big(\Omega;[C^1(\overline{D_\refd})]^3\big)
\subset L^2\big(\Omega;[C^1(\overline{D_\refd})]^3\big)\), 
the deformation field \({\bs\chi}_{D}\) can be represenetd
by a Karhunen-Lo\`eve expansion \cite{Loe77} which has the 
form
\begin{equation}\label{eq:randomVectorFieldDomain}
{\bs\chi}_{D}(\xref,\omega)=\E[{\bs\chi}_{D}](\xref)+
\sum_{k=1}^\infty\sqrt{\lambda_{D,k}}{\bs\chi}_{D,k}(\xref){Y}_{D,k}(\omega),\quad
\xref\in D_\refd.
\end{equation}
Herein,
\[
\E[{\bs\chi}_{D}](\xref)\isdef\int_\Omega {\bs\chi}_{D}(\xref,\omega)\d\mathbb{P}(\omega)
\]
denotes the expectation, while \((\lambda_{D,k}, {\bs\chi}_{D,k})\) are 
the eigenpairs of the covariance operator $\mathcal{C}_{D}\colon \big[L^2(D)\big]^3\to \big[L^2(D)\big]^3$,
\begin{align}\label{eq:domaincovarianceoperator}
(\mathcal{C}_{D}{\bs U})(\xref)\isdef\int_{{D_\refd}} \Cov[{\bs \chi}_{D}](\xref,\xref'){\bs U}(\xref')\d\xref',
\end{align}
where
\[
\Cov[{\bs\chi}_{D}](\xref,\xref')\isdef\int_\Omega 
\big({\bs\chi}_{D}(\xref,\omega)-\E[{\bs\chi}_{D}](\xref)\big)
\big({\bs\chi}_{D}(\xref',\omega)-\E[{\bs\chi}_{D}](\xref')\big)^\intercal\d\mathbb{P}(\omega).
\]
It holds that
\begin{equation}\label{eq:randvar}
Y_{D,k}(\omega)\isdef\frac{1}{\sqrt{\lambda_{D,k}}}\int_{D_\refd}
\big({\bs\chi}_{D}(\xref,{\omega})-\E[{\bs\chi}_{D}](\xref)\big)^\intercal{\bs\chi}_{D,k}(\xref)\d\xref.
\end{equation}
The family \(\{Y_{D,k}\}_k\)  of random variables is therefore uncorrelated and centred. We remark that in uncertainty quantification problems typically only \(\E[{\bs\chi}_{D}]\) and \(\Cov[{\bs \chi}_{D}]\)
are known, such that the random variables cannot be inferred via \eqref{eq:randvar}. Instead, their (common) distribution has to be appropriately estimated.

\subsection{Modelling of random surfaces}
The numerical computation of a Karhunen-Lo\`eve expansion as outlined in the previous subsection will generally require a (volume) finite element mesh for \(D_\refd\). Moreover, the data \(\E[{\bs\chi_{D}}]\) and \(\Cov[{\bs \chi_{D}}]\) need to be known on the whole reference domain $D_\refd$. In contrast, for our boundary element-based approach, we only require realizations of the perturbed boundary. Especially, the following exposition shows that, for the 
computation of surface realizations, knowledge of 
\(\E[{\bs\chi}_{D}]\) and \(\Cov[{\bs \chi}_{D}]\) at the boundary $S_\refd=\partial D_\refd$ is sufficient.

Given a function 
\(\mathbf{g}\colon D_\refd\to\mathbb{R}^3\), let
\[
(\Dtrace \mathbf{g})(\xref)\isdef\lim_{\xref'\ni 
D_\refd\to \xref\in S_\refd}\mathbf{g}(\xref')
\]
denote the (interior) trace operator and ${\bs\chi}_{S}\isdef\Dtrace{\bs\chi}_{D}$. Since
\[
\Dtrace\colon\big[C^1\big(\overline{D_\refd}\big)\big]^3\subset \big[H^1(D_\refd)\big]^3\to \big[H^{1/2}(S_\refd)\big]^3
\]
is a
continuous operator and the Bochner integral commutes 
with continuous operators, \cite{HP57}, it holds
\[
(\Dtrace\E[{\bs \chi}_{D}])(\xref)
 = \Dtrace\int_\Omega {\bs\chi}_{D}(\xref,\omega)\d\mathbb{P}(\omega)
 = \int_\Omega {\bs\chi}_{S}(\xref,\omega)\d\mathbb{P}(\omega)
 = \E[{\bs \chi}_S](\xref)
\]
as well as
\begin{align*}
&(\Dtrace\otimes \Dtrace)\Cov[{\bs\chi}_{D}](\xref,\xref')\\
&\quad=\int_\Omega \big((\Dtrace{\bs\chi}_{D})(\xref,\omega)-(\Dtrace\E[{\bs\chi}_{D}])(\xref)\big)\\
&\hspace{6cm}\cdot\big((\Dtrace{\bs\chi}_{D})(\xref',\omega)-
(\Dtrace\E[{\bs\chi}_{D}])(\xref')\big)^\intercal\d\mathbb{P}(\omega)\\
&\quad=\int_\Omega \big({\bs\chi}_{S}(\xref,\omega)-(\E[{\bs\chi}_{S}])(\xref)\big)
\big({\bs\chi}_{S}(\xref',\omega)-
\E[{\bs\chi}_{S}](\xref')\big)^\intercal\d\mathbb{P}(\omega)\\
&\quad=\Cov[{\bs\chi}_{S}](\xref,\xref').
\end{align*}
Therefore, the random 
deformation field at \(S\), i.e.\ \({\bs\chi}_{S}
(\xref,\omega)\), is fully described by \((\Dtrace
\E[{\bs \chi}_{D}])(D)\) and \((\Dtrace\otimes
\Dtrace)\Cov[{\bs\chi}_{D}](\xref,\xref')\). For the numerical computation of the deformation field, it is therefore sufficient to compute the eigenpairs \((\lambda_{S,k}, {\bs\chi}_{S,k})\) of the surface covariance operator $\mathcal{C}_S\colon \big[L^2(S_\refd)\big]^3\to \big[L^2(S_\refd)\big]^3$ given by
\begin{align}\label{eq:surfacecovarianceoperator}
(\mathcal{C}_{S_\refd}{\bs U})(\xref)\isdef\int_{S_\refd} (\Dtrace\otimes \Dtrace)\Cov[{\bs \chi}_{D}](\xref,\xref'){\bs U}(\xref')\d\sigma_{\xref'},
\end{align}
to obtain 
\begin{equation}\label{eq:randomVectorFieldBoundary}
{\bs\chi}_{S}(\xref,\omega)=\Dtrace\E[{\bs\chi}_{D}](\xref)+
\sum_{k=1}^\infty\sqrt{\lambda_{S,k}}{\bs\chi}_{S,k}(\xref){Y}_{S,k}(\omega)
\end{equation}
with
\[
Y_{S,k}(\omega)\isdef\frac{1}{\sqrt{\lambda_{S,k}}}\int_{S}
\big({\bs\chi}_{S}(\xref,{\omega})-\Dtrace\E[{\bs\chi}_{D}](\xref)\big)^\intercal{\bs\chi}_{S,k}(\xref)\d\xref.
\]

We remark that the computation of the eigenpairs of \eqref{eq:surfacecovarianceoperator} is significantly cheaper than the computation of the ones of \eqref{eq:domaincovarianceoperator} since the latter only relies on a surface mesh for $S_{\refd}$ rather than a volume mesh for $D_{\refd}$. Thus, the discrete system will be significantly smaller.
However, that the corresponding eigenfunctions will in general not be 
traces of the eigenfunctions of the Karhunen-Lo\`eve expansion 
\eqref{eq:randomVectorFieldDomain} and also the distribution of the random 
variables will change. In the sequel, we assume that the family
\(\{Y_{S,k}\}_k\) is independent and uniformly distributed with
\(\{Y_{S,k}\}_k\sim\mathcal{U}(-1,1)\) for all \(k\). Then, we can
identify each of the random variables by its image \(y_k\in[-1,1]\)
and up with the parametric deformation field
\[
{\bs\chi}_{S}(\xref,{\bs y})=\Dtrace\E[{\bs\chi}_{D}](\xref)+
\sum_{k=1}^\infty\sqrt{\lambda_{S,k}}{\bs\chi}_{S,k}(\xref)y_k,\quad
{\bs y}\in\Gamma\isdef[-1,1]^{\mathbb{N}},
\]
which gives rise to the random surface
\begin{align}\label{eq:randomdomain}
  S({\bs y}) = \bigg\{{\bs\chi}_S(\xref,{\bs y}):
\xref\in S_\refd\bigg\}.
\end{align}

\section{Isogeometric discretization of random domains}\label{sec:DiscRD}

\subsection{Fundamental Notions}\label{sec::subsec::iga}
We review the basic notions of isogeometric analysis, restricting ourselves to spaces constructed via locally quasi-uniform $p$-open knot vectors as required by the theory presented in \cite{Beirao-da-Veiga_2014aa,BDK+2020}.

\begin{definition}\label{def::splines}
	Let $p$ and $k$ with $0\leq p< k$. 
	A \emph{locally quasi uniform $p$-open knot vector} is a tuple
	\begin{align*}
	\Xi = \big[{\xi_0 = \cdots =\xi_{p}}\leq \cdots \leq{\xi_{k}=\cdots =\xi_{k+p}}\big]\in[0,1]^{k+p+1}
	\end{align*}
	with $\xi_0 = 0$ and $\xi_{k+p}=1$
	such that there exists a constant $\theta\geq 1$ with
	 $\theta^{-1}\leq h_j\cdot h_{j+1}^{-1} \leq \theta$
	 for all $p\leq j < k$, 
	  where  $h_j\isdef  \xi_{j+1}-\xi_{j}$.
	The B-spline basis $ \lbrace b_j^p \rbrace_{0\leq j< k}$ is then recursively
	defined according to
	\begin{align*}
	b_j^p(x) & =\begin{cases}
	\mathbbm{1}_{[\xi_j,\xi_{j+1})}&\text{ if }p=0,\\[8pt]
	\frac{x-\xi_j}{\xi_{j+p}-\xi_j}b_j^{p-1}(x) +\frac{\xi_{j+p+1}-x}{\xi_{j+p+1}-\xi_{j+1}}b_{j+1}^{p-1}(x) & \text{ else,}
	\end{cases}
	\end{align*}
	where $\mathbbm{1}_A$ refers to the indicator function of the set $A$. The corresponding
	spline space is finally defined according to $\S^p(\Xi)\isdef\operatorname{span}(\lbrace b_j^p\rbrace_{j <k}).$
\end{definition}

To obtain spline spaces in two spatial dimensions, we employ a tensor product construction.
More precisely, for a tuple $\bfXi =(\Xi_1,$ $\Xi_2)$ and polynomial degrees $\bfp=(p_1,p_2)$, we define the spaces \[
\S^{\bfp}(\bfXi)\isdef \S^{p_1}(\Xi_1)\otimes \S^{p_2}(\Xi_2).
\]
Given knot vectors $\Xi_1,$ $\Xi_2$ with knots $\xi_{i}^k < \xi_{i+1}^k$ for $k=1,2$, sets of the form $[\xi_{j}^1,\xi^1_{j+1}]\times[\xi^2_{j},\xi^2_{j+1}]$ will be called \emph{elements}. We reserve the letter $h$ for the maximal diameter of all elements.
For further concepts and algorithmic realization of B-splines, we refer to \cite{Piegl_1997aa} and the references therein.

\subsection{Isogeometric boundary representation}\label{sec:surfrep}
In the following, we will assume the usual isogeometric setting for the surface
$S_\refd$ of the reference domain $D_\refd$, i.e.\ denoting the unit square by 
$\square\isdef[0,1]^2$, we assume that the surface $S_\refd$ can be decomposed 
into several smooth \emph{patches}
\[
S_\refd = \bigcup_{i=1}^M S_{\refd}^{(i)},
\]
where the intersection $S_{\refd}^{(i)}\cap S_{\refd}^{(i')}$ 
consists at most of a common vertex or a common edge for 
\(i\neq i^\prime\). In particular, we model each patch $S_{\refd}^{(i)}$ as 
an invertible NURBS mapping
\begin{equation}\label{eq:parametrization}
{\bs s}_i\colon\square\to S_{\refd}^{(i)}\quad\text{ with }\quad S_{\refd}^{(i)} = {\bs s}_i(\square)
\quad\text{ for } i = 1,2,\ldots,M,
\end{equation}
where \({\bs s}_i\) is of the form
\begin{align*}
\mathbf{s}_i(x,y)\isdef \sum_{0=i_1}^{k_1}\sum_{0=i_2}^{k_2}\frac{\mathbf{c}_{i_1,i_2} b_{i_1}^{p_1}(x) b_{i_2}^{p_2}(y) w_{i_1,i_2}}{ \sum_{j_1=0}^{k_1-1}\sum_{j_2=0}^{k_2-1} b_{j_1}^{p_1}(x) b_{j_2}^{p_2}(y) w_{j_1,j_2}}
\end{align*}
for control points $\mathbf{c}_{i_1,i_2}\in \Rset^3$ and weights $w_{i_1,i_2}>0$. 
We shall further follow the common convention that
parametrizations with a common edge coincide except for orientation.
\begin{figure}[htb]
	\begin{center}
		\begin{tikzpicture}[
		scale=.3,
		axis/.style={thick, ->, >=stealth'},
		important line/.style={thick},
		every node/.style={color=black}
		]
		\draw (-.5,-1)node{{$0$}};
		\draw (8,-1)node{{$1$}};
		\draw (-.5,8)node{{$1$}};
		\draw (0,0)--(8,0);
		\draw (0,1)--(8,1);
		\draw (0,2)--(8,2);
		\draw (0,3)--(8,3);
		\draw (0,4)--(8,4);
		\draw (0,5)--(8,5);
		\draw (0,6)--(8,6);
		\draw (0,7)--(8,7);
		\draw (0,8)--(8,8);
		\draw (0,0)--(0,8);
		\draw (1,0)--(1,8);
		\draw (2,0)--(2,8);
		\draw (3,0)--(3,8);
		\draw (4,0)--(4,8);
		\draw (5,0)--(5,8);
		\draw (6,0)--(6,8);
		\draw (7,0)--(7,8);
		\draw (8,0)--(8,8);
		\draw (3,3)node(N1){};
		\draw (19.3,7)node(N2){};
		\draw
		(22.5,4)node{\includegraphics[width=0.3\textwidth,clip=true,trim=185 540
			180 90]{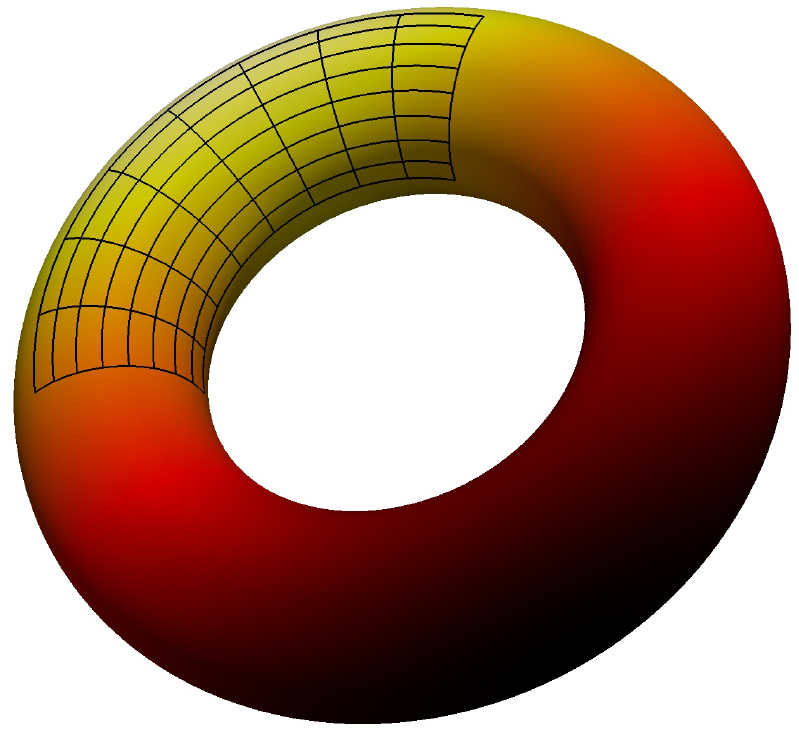}};
		\draw (11.5,8.6)node(N3){\large${\bs s}_i$};
		\draw (21,5)node(N4){{\large${S}_\refd^{(i)}$}};
		\path[->,line width=1pt]  (N1) edge [bend left] (N2);
		\end{tikzpicture}
	\end{center}
	\caption{Surface representation and mesh generation.}
	\label{fig:mesh}
\end{figure}
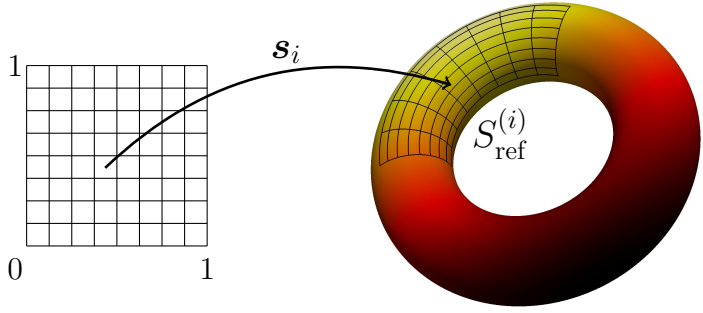

Following the spirit of isogeometric analysis, the random surface $S(\bs{y})$ from \eqref{eq:randomdomain} will be represented as a union of NURBS patches. This
is achieved by appropriately discretizing the deformation field \eqref{eq:randomVectorFieldBoundary}. More precisely, the random surface
$S(\bs{y})$ is discretized by $S(\bs{y})\approx S_h(\bs{y})$,
where the latter can be decomposed into $M$ distinct NURBS patches
\[
S_h(\bs{y}) = \bigcup_{i=1}^M S_h^{(i)}(\bs{y}).
\]
Herein, the intersection $S_h^{(i)}(\bs{y})\cap S_h^{(i)}(\bs{y})$ again
consists at most of a common vertex or a common edge for 
\(i\neq i^\prime\) and each patch $S_h^{(i)}(\bs{y})$ is given by an invertible mapping
\[
{\bs s}_{i,h}(\cdot,\bs{y})\colon\square\to S^{(i)}(\bs{y})\quad\text{ with }\quad S^{(i)}(\bs{y}) = {\bs s}_{i,h}(\square,\bs{y})
\quad\text{ for } i = 1,2,\ldots,M,
\]
with
\[
\bs{s}_{i,h}(\xref,\bs{y}) = \bs{s}_i(\xref)+{\bs\chi}_{S,h}\big|_i(\xref,\bs{y}).
\]
First, we note that $\mathbf{s}_{i,h}(\xref,\bs{y})$ is again a NURBS mapping if ${\bs\chi}_{S,h}\big|_i(\xref,\mathbf{y})$ is discretized by using appropriate basis functions. In fact, if these basis functions are chosen also as NURBS, the randomness of the surface is transformed onto transformations of the control points. Second, we note that ${\bs\chi}_{S}$ needs to be at least globally continuous to obtain an admissible surface transformation. Given a tuple of knot vectors $\bfXi$ and polynomial degrees $\bfp$, a natural choice for the discretization of ${\bs\chi}_{S}(\cdot, \bs{y})$ is thus given by the vector valued spline space
\[
\bbS_{\bfp,\bfXi}(S_\refd)=\big[\S_{\bfp,\bfXi}(S_\refd)\big]^3
\]
where
\[
\S_{\bfp,\bfXi}(S_\refd)\isdef \left\lbrace f\in C(S_\refd)\colon f|_i\circ \bs{s}_i^{-1}\in \S_{\bfp}(\bfXi)\text{ for }1\leq i\leq M\right\rbrace.
\]
Of course, the knot vectors and polynomial degrees could vary in each component and on each patch, but for simplicity we choose to use the same knots and degrees for better readability. Approximation properties for these spaces where derived in \cite{BDK+2020}.

Next, we discuss how an approximation of ${\bs\chi}_{S,h}$ in terms of such basis functions can be derived
by computing Karhunen-Lo\`eve expansion \eqref{eq:randomVectorFieldBoundary} of the underlying random deformation.

\subsection{Fast computation of the Karhunen-Lo\`eve expansion}

The computation of the Karhunen-Lo\`eve expansion of surface deformations from
the expectation and the covariance is directly related to the solution of the 
eigenvalue problem
\[
\mathcal{C}_{S}{\bs\chi}_{S,k}=\lambda_{S,k}{\bs\chi}_{S,k}.
\]
Based on the previous discussion, 
it is natural to choose a B-spline-based Galerkin discretization for the numerical
solution of the eigenvalue problem. Hence, replacing $\big[L^2(S_\refd)\big]^3$ by a B-spline space $\bbS_{\bfp,\bfXi}(S_\refd)$ in the eigenproblem's weak formulation
\begin{align*}
&\text{Find $(\lambda_{S,k},{\bs\chi}_{S,k})\in\mathbb{R}\times \big[L^2(S_\refd)\big]^3$ such that}\\
&\quad(\mathcal{C}_S{\bs\chi}_{S,k},v)_{[L^2(S_\refd)]^3}=\lambda_{S,k}({\bs\chi}_{S,k},v)_{[L^2(S_\refd)]^3}
\quad\text{for all $v\in \big[L^2(S_\refd)\big]^3$}
\end{align*}
yields the discrete generalized eigenvalue problem
\begin{equation}\label{eq:disceig}
\underline{\bfC}\underline{\boldsymbol{\chi}}_{k}=\lambda_{k,h}\underline{\bfM}\underline{\boldsymbol{\chi}}_{k}.
\end{equation}
Although the mass matrix $\underline{\mathbf{M}}$ is sparse, the covariance matrix $\underline{\bfC}$ is typically densely populated,
as it issues from the discretization of a nonlocal operator. Therefore, a naive solution 
of this eigenvalue problem is prohibitive for a larger number of degrees of freedom.

As a viable alternative, we assume that a low-rank factorization $\underline{\bfC}\approx\underline{\bfL}\underline{\bfL}^\intercal$ of the covariance matrix is known. Such a factorization can, for example, efficiently be computed by the truncated pivoted Cholesky decomposition, see \cite{HPS}. Inserting this decomposition into \eqref{eq:disceig} yields
\begin{equation}\label{eq:bigpceig}
\underline{\bfL}\underline{\bfL}^{\intercal}\underline{\boldsymbol{\chi}}_{k}=\lambda_{k,h}\underline{\bfM}\underline{\boldsymbol{\chi}}_{k}.
\end{equation}
Substituting $\underline{\bfpsi}_k=\underline{\mathbf{M}}^{1/2}\underline{\bfchi}_k$ therefore results in the eigenvalue problem
\begin{equation}\label{eq:smallpceig}
\underline{\bfL}^{\intercal}\underline{\bfM}^{-1}\underline{\bfL}\underline{\bfpsi}_{k}=\lambda_{k,h}\underline{\bfpsi}_{k},
\end{equation}
which has the same non-zero eigenvalues as \eqref{eq:bigpceig}, but is significantly smaller and cheaper to compute if $\underline{\bfC}$ has low rank. The eigenvectors of \eqref{eq:bigpceig} can be retrieved from \eqref{eq:smallpceig} by making use of the relation $\bfchi_k=\underline{\mathbf{M}}^{-1}\underline{\bfL}\bfpsi_k.$
\begin{remark}\label{rem:efficientIGAKL}
The supports of the basis functions in $\bbS_{\bfp,\bfXi}(S_\refd)$ can be quite large, which makes the assembly of a single matrix entry as used for the truncated pivoted Cholesky decomposition computationally expensive. Instead, one may opt for performing the Cholesky decomposition directly on the matrix $\underline{\bfC}_\star$ generated by the shape functions $\bbS_{\bfp,\bfXi}^\star(S_\refd)$ of $\bbS_{\bfp,\bfXi}(S_\refd)$. Then, there exists a matrix version $\underline{\bfT}$ of the local-to-global map such that $\underline{\bfC}=\underline{\bfT}\underline{\bfC}_\star\underline{\bfT}^\intercal$. Now, substituting $\underline{\bfC}_\star\approx\underline{\bfL}_\star\underline{\bfL}_\star^\intercal$ yields a low-rank factorization $\underline{\bfC}\approx\underline{\bfT}\underline{\bfL}_\star(\underline{\bfT}\underline{\bfL}_\star)^\intercal=\underline{\bfL}\underline{\bfL}^\intercal$.
\end{remark}

\section{Boundary integral equations}\label{sec:bie}
\subsection{Computing the scattered wave}
We recall the solution of the boundary value problem 
\eqref{eq:pde} by means of boundary integral equations. 
To this end, and for sake of simplicity in representation, we 
assume for the moment that the domain $D$ is fixed and has a 
Lipschitz surface $S = \partial D$.

We introduce the acoustic single layer operator
\[
  \mathcal{V}\colon H^{-\nicefrac{1}{2}}(S)\to H^{\nicefrac{1}{2}}(S),\quad
  \mathcal{V}\rho\isdef\int_S
  	\Phi(\cdot,{\bs z})\rho({\bs z})\d\sigma_{\bs z}
\]
and the acoustic double layer operator
\[
  \mathcal{K}\colon L^2(S)\to L^2(S),\quad
  \mathcal{K}\rho\isdef\int_S
  	\frac{\partial\Phi(\cdot,{\bs z})}{\partial{\bs n}_{\bs z}}\rho({\bs z})\d\sigma_{\bs z}.
\]
Here, ${\bs n}_{\bs z}$ denotes the outward pointing normal vector
at the surface point ${\bs z}\in S$,  while $\Phi(\cdot,\cdot)$ denotes 
the Green's function for the Helmholtz equation. In three spatial 
dimensions, the Green's function is given by
\[
  \Phi({\bs x}, {\bs z}) = \frac{e^{i\kappa\|{\bs x}-{\bs y}\|_2}}{4\pi\|{\bs x}-{\bs y}\|_2}.
\]
Considering an incident plane wave $u_{\mathrm{inc}}({\bs x}) = e^{i\kappa\langle{\bs d},
{\bs x}\rangle}$, $\|{\bs d}\|_2=1$, the Neumann data of the total wave $u=u_{\mathrm{inc}}+u_{\mathrm{s}}$ at the 
surface $S$ can be determined by the boundary integral equation
\begin{equation}\label{eq:BIE1}
  \left(\frac{1}{2} + \mathcal{K^\star} - i\eta\mathcal{V}\right)
  	\frac{\partial u}{\partial {\bs n}} = \frac{\partial u_{\mathrm{inc}}}{\partial {\bs n}} - i\eta u_{\mathrm{inc}}
		\quad\text{on $S$},
\end{equation}
with $\eta=\kappa/2$, compare to \cite{CK2}. 

From the Cauchy data of $u$ at $S$, we 
can determine the scattered wave $u_{\mathrm{s}}$ in any point in the
exterior of the obstacle by applying the potential evaluation
\begin{equation}\label{eq:solution1}
  u_{\mathrm{s}}({\bs x}) = \int_S\Phi({\bs x}, {\bs z})
	\frac{\partial u}{\partial{\bs n}}({\bs z})\d\sigma_{\bs z},
	\quad {\bs x}\in\mathbb{R}^3\setminus\overline{D}.
\end{equation}

\subsection{Scattered wave representation at an artificial interface}
We introduce an artificial interface $T\subset\mathbb{R}^3$, 
being sufficiently large to guarantee that $T$ encloses all realizations of
the domain $D$. In view of \eqref{eq:solution1}, we may compute the 
Cauchy data $u_{\mathrm{s}}|_T$ and $(\partial u_{\mathrm{s}}/\partial {\bs n})|_T$ 
of the scattered wave at the artificial interface $T$. It holds
\[
  \frac{\partial u_{\mathrm{s}}}{\partial {\bs n}}({\bs x}) 
  	= \int_S \frac{\partial\Phi({\bs x}, {\bs z})}
  		{\partial {\bs n}_{\bs x}}\frac{\partial u}{\partial{\bs n}}({\bs z})\d\sigma_{\bs z},
				\quad {\bs x}\in T.
\]

For any ${\bs x}\in\mathbb{R}^3$ located outside the artificial 
interface, we may now either employ the representation formula 
\eqref{eq:solution1} or the representation formula
\begin{equation}\label{eq:solution2}
  u_{\mathrm{s}}({\bs x}) = \int_T \bigg\{\Phi({\bs x}, {\bs z})
	\frac{\partial u_{\mathrm{s}}}{\partial{\bs n}}({\bs z})+
	\frac{\partial\Phi({\bs x}, {\bs z})}{\partial {\bs n}_{\bs z}}
		u_{\mathrm{s}}({\bs z})\bigg\}\d\sigma_{\bs z}
\end{equation}
to compute the scattered wave $u_{\mathrm{s}}$.

The major advantage of \eqref{eq:solution2} over \eqref{eq:solution1}
is that the artificial interface is fixed in contrast to the shape of the 
random obstacle later on. 

\subsection{Scattering at random obstacles}\label{sec:random}
From now on, let the obstacle be subject to uncertainty as introduced in
Section~\ref {sec:randomDomains}. We describe the uncertain
obstacle $D = D({\bs y})$ by its {random surface} $S({\bs y})$,
which is given by \eqref{eq:randomdomain}.

Having the incident wave $u_{\mathrm{inc}}$ at hand, the boundary 
value problem for the total field $u({\bs y}) = u_{\mathrm{s}}({\bs y})+u_{\mathrm{inc}}$
for any ${\bs y}\in\Gamma$
reads 
\begin{equation}\label{eq:spde}
 \begin{aligned}
  \Delta u({\bs y}) + \kappa^2 u({\bs y}) = 0\quad&\text{in}\ \mathbb{R}^3\setminus\overline{D({\bs y})},\\
  	u({\bs y}) = 0\quad&\text{on}\ S({\bs y}),\\
	  \sqrt{r}\bigg(\frac{\partial u_{\mathrm{s}}}{\partial r}-i\kappa u_{\mathrm{s}}\bigg)
  	\to 0\quad&\text{as}\ r = \|{\bs x}\|_2\to\infty.
  \end{aligned}
\end{equation}
By the construction of \(S({\bs y})\), the random scattering 
problem \eqref{eq:spde} exhibits a unique solution for each realization 
\({\bs y}\in\Gamma\) of the random parameter.
Moreover, it has been shown in \cite{H3S} for the case of the Helmholtz transmission 
problem that the total wave \(u({\bs y})\) exhibits an analytic extension into a 
certain region of the complex plane with respect to the parameter 
\({\bs y}\in\Gamma\). This particularly
allows for the use of higher order quadrature methods, like quasi-Monte 
Carlo methods, see e.g.\ \cite{Caf98,NIE}, or even sparse quadrature methods, see 
e.g.\ \cite{HHPS18,H3S} in order to compute quantities of interest, such as 
expectation and variance. Extensions to the Maxwell case are discussed in \cite{JHS17,JSZ17,AJZ20}.

\subsection{Expectation of the scattered wave}
The scattered wave's expectation can be computed 
for any given point ${\bs x}\in\mathbb{R}^3$ by the 
representation formula \eqref{eq:solution1}, 
which leads to
\begin{equation}\label{eq:expectation1}
  \E[u_{\mathrm{s}}]({\bs x}) = \E\bigg[\int_{S({\bs y})}\Phi({\bs x}, {\bs z})
	\frac{\partial u_{\mathrm{s}}}{\partial{\bs n}}({\bs z},\cdot)\d\sigma_{\bs z}\bigg].
\end{equation}
Obviously, \eqref{eq:expectation1} only makes sense if ${\bs x}
\in\mathbb{R}^3$ is sufficiently far away from the random obstacle.
Otherwise, there might be instances ${\bs y}\in\Gamma$ such that ${\bs x}\in D({\bs y})$, i.e.\ 
the point ${\bs x}\in\mathbb{R}^3$ does not lie outside the 
obstacle almost surely.

If the expectation needs to be evaluated at many locations, it is 
much more efficient to introduce the artificial interface 
$T$ and to consider expression \eqref{eq:solution2}. 
For any ${\bs x}\in\mathbb{R}^3$ lying outside the interface 
$T$, it holds
\begin{equation}\label{eq:expectation2}
  \E[u_{\mathrm{s}}]({\bs x}) = \int_T \bigg\{\Phi({\bs x}, {\bs z})
	\E\bigg[\frac{\partial u_{\mathrm{s}}}{\partial{\bs n}}\bigg]({\bs z})+
	\frac{\partial\Phi({\bs x}, {\bs z})}{\partial {\bs n}_{\bs z}}
		\E[u_{\mathrm{s}}]({\bs z})\bigg\}\d\sigma_{\bs z}.
\end{equation}
As a consequence, the scattered wave's expectation is completely
encoded in the Cauchy data at the artificial interface $T$.
This means that we only need to compute the expected Cauchy data 
\begin{equation}\label{eq:expectation3}
  \mathbb{E}[u_{\mathrm{s}}] = \int_{\Gamma} 
  	\bigg\{\int_{S({\bs y})}\Phi({\bs x}, {\bs z})
	\frac{\partial u}{\partial{\bs n}}({\bs z},{\bs y})\d\sigma_{\bs z}\bigg\}\d\mu
\end{equation}
and
\begin{equation}\label{eq:expectation4}
  \mathbb{E}\bigg[\frac{\partial u_{\mathrm{s}}}{\partial {\bs n}}\bigg]
  	= \int_{\Gamma} 
  		\bigg\{\int_{S({\bs y})}\frac{\Phi({\bs x}, {\bs z})}{\partial {\bs n}_{\bs z}}
			u({\bs z},{\bs y})\d\sigma_{\bs z}\bigg\}\d\mu
\end{equation}
of the scattered wave at the artificial interface $T$, 
which is of lower spatial dimension than the exterior domain.

\subsection{Variance of the scattered wave}
The variance $\V[u_{\mathrm{s}}]$ of the scattered wave $u_{\mathrm{s}}({\bs y)}$ 
at a point ${\bs x}\in\mathbb{R}^3$ outside the artificial interface 
$T$ depends nonlinearly on the Cauchy data of $u_{\mathrm{s}}$ 
at the interface. Nonetheless, we can make use of the fact that the 
variance is the trace --in the algebraic sense-- of the covariance, i.e.\
\begin{equation}\label{eq:var}
  \V[u_{\mathrm{s}}]({\bs x}) = \Cov[u_{\mathrm{s}}]({\bs x},{\bs x}')\big|_{{\bs x}={\bs x}'}\\
  	= \Cor[u_{\mathrm{s}}]({\bs x},{\bs x}')\big|_{{\bs x}={\bs x}'} - |\E[u_{\mathrm{s}}]({\bs x})|^2,
\end{equation}
where the covariance is given by
\begin{align*}
  \Cov[u_{\mathrm{s}}]({\bs x},{\bs x}') &= \E\Big[\big(u_{\mathrm{s}}({\bs x},\cdot)-\E[u_{\mathrm{s}}]({\bs x})\big)
  	\overline{\big(u_{\mathrm{s}}({\bs x}',\cdot)-\E[u_{\mathrm{s}}]({\bs x}')\big)}\Big]\\
	&=\E\big[u_{\mathrm{s}}({\bs x},\cdot)\overline{u_{\mathrm{s}}({\bs x'},\cdot)}\big]-\E[u_{\mathrm{s}}]({\bs x})\overline{\E[u_{\mathrm{s}}]({\bs x}')}.
\end{align*}
Hence, it holds for the correlation
\[\Cor[u_{\mathrm{s}}]({\bs x},{\bs x}') = \E\big[u_{\mathrm{s}}({\bs x},\cdot)\overline{u_{\mathrm{s}}({\bs x'},\cdot)}\big].
\]
The correlation is a higher-dimensional 
object which depends only linearly on the second 
moment of the Cauchy data of the scattered wave at the 
artificial interface $T$. This greatly simplifies the computation
of the variance. More precisely, by defining for ${\bs x},{\bs x}'\in T$ 
the correlations
\begin{align*}
  \Cor[u_{\mathrm{s}}]({\bs x},{\bs x}')
  &= \E\bigg[\bigg(\int_{S({\bs y})}\!\!\!\!\!\Phi({\bs x}, {\bs z})
	\frac{\partial u_{\mathrm{s}}}{\partial{\bs n}}({\bs z},{\bs y})\d\sigma_{\bs z}\bigg)
	\overline{\bigg(\int_{S({\bs y})}\!\!\!\!\!\Phi({\bs x}', {\bs z})
	\frac{\partial u_{\mathrm{s}}}{\partial{\bs n}}({\bs z},{\bs y})\d\sigma_{{\bs z}}\bigg)}\bigg],\\
  \Cor\bigg[\frac{\partial u_{\mathrm{s}}}{\partial{\bs n}}\bigg]({\bs x},{\bs x}')
  &= \E\bigg[\bigg(\int_{S({\bs y})}\!\!\!\!\!\frac{\partial\Phi({\bs x}, {\bs z})}{\partial {\bs n}_{\bs z}}
		u_{\mathrm{s}}({\bs z},{\bs y})\d\sigma_{\bs z}\bigg)
	\overline{\bigg(\int_{S({\bs y})}\!\!\!\!\!\frac{\partial\Phi({\bs x}', {\bs z})}{\partial {\bs n}_{\bs z}}
		u_{\mathrm{s}}({\bs z},{\bs y})\d\sigma_{{\bs z}}\bigg)}\bigg],
\end{align*}
and
\begin{align*}
  &\Cor\bigg[u_{\mathrm{s}},\frac{\partial u_{\mathrm{s}}}{\partial{\bs n}}\bigg]({\bs x},{\bs x}')
  	= \overline{\Cor\bigg[\frac{\partial u_{\mathrm{s}}}{\partial{\bs n}},u_{\mathrm{s}}\bigg]({\bs x}',{\bs x})}\\
  &\qquad= \E\bigg[\overline{\bigg(\int_{S({\bs y})}\Phi({\bs x}, {\bs z})
	\frac{\partial u_{\mathrm{s}}}{\partial{\bs n}}({\bs z},\omega)\d\sigma_{\bs z}\bigg)}
	\bigg(\int_{S({\bs y})}\frac{\partial\Phi({\bs x}', {\bs z})}{\partial {\bs n}_{\bs z}}
		u_{\mathrm{s}}({\bs z},{\bs y})\d\sigma_{{\bs z}}\bigg)\bigg],
\end{align*}
we find for two points ${\bs x},{\bs x}'\in\mathbb{R}^3$ lying 
outside of the interface $T$ the deterministic expression 
\begin{equation}\label{eq:cor}
\begin{aligned}
  \Cor[u_{\mathrm{s}}]({\bs x},{\bs x}') = \int_T\int_T
  	\bigg\{&\Phi({\bs x}, {\bs z})\overline{\Phi({\bs x}', {\bs z}')}\Cor\!\!\bigg[\frac{\partial u_{\mathrm{s}}}{\partial{\bs n}}\bigg]({\bs z},{\bs z}')\\
  &\ + \Phi({\bs x}, {\bs z})\overline{\frac{\partial\Phi({\bs x}', {\bs z}')}{\partial {\bs n}_{{\bs z}'}}}
  	\Cor\!\!\bigg[\frac{\partial u_{\mathrm{s}}}{\partial{\bs n}},u_{\mathrm{s}}\bigg]({\bs z},{\bs z}')\\
  &\ + \frac{\partial\Phi({\bs x}, {\bs z})}{\partial {\bs n}_{\bs z}}\overline{\Phi({\bs x}', {\bs z}')}
  	\Cor\!\!\bigg[u_{\mathrm{s}},\frac{\partial u_{\mathrm{s}}}{\partial{\bs n}}\bigg]({\bs z},{\bs z}')\\
  &\ + \frac{\partial\Phi({\bs x}, {\bs z})}{\partial {\bs n}_{\bs z}}
	\overline{\frac{\partial\Phi({\bs x}', {\bs z}')}{\partial {\bs n}_{{\bs z}'}}}
	\Cor[u_{\mathrm{s}}]({\bs z},{\bs z}')\bigg\}\d\sigma_{{\bs z}'}\d\sigma_{\bs z}.
\end{aligned}
\end{equation}

\section{Multilevel quadrature}\label{sec:MLQ}
In order to calculate quantities of interest efficiently, we employ a multilevel quadrature approach. For the computation of the expectation,
we may exploit the linearity of the expectation in formula \eqref{eq:expectation2} and rely
on the Cauchy data on the spatial refinement
levels \(\ell=0,1,\ldots,L\) computed at the artificial interface \(T\).
Thus, we obtain
\begin{equation}\label{eq:MLQmean}
  \E[u_{\mathrm{s}}]({\bs x}) 
		\approx \int_T \bigg\{\Phi({\bs x}, {\bs z})
	\mathcal{Q}_{L}^{\text{ML}}\bigg[\frac{\partial u_{\mathrm{s}}}{\partial{\bs n}}\bigg]({\bs z})+
	\frac{\partial\Phi({\bs x}, {\bs z})}{\partial {\bs n}_{\bs z}}
		\mathcal{Q}_{L}^{\text{ML}}[u_{\mathrm{s}}]({\bs z})\bigg\}\d\sigma_{\bs z}
\end{equation}
with
\[
\mathcal{Q}_{L}^{\text{ML}}[\rho]({\bs z})\isdef\sum_{\ell=0}^L
\mathcal{Q}_{L-\ell}\big(\rho^{(\ell)}({\bs z},\cdot)-
\rho^{(\ell-1)}({\bs z},\cdot)\big)\quad\text{for }{\bs z}\in T,
\]
where \(\mathcal{Q}_{\ell}\) is a quadrature rule on level \(\ell\).
Moreover, the function \(\rho^{(\ell)}\) is the Galerkin projection of the
density \(\rho\) evaluated at the artificial interface for the spatial
refinement on level \(\ell\) of the scatterer, where we set \(\rho^{(-1)}\equiv 0\).

For the approximation error of the multilevel quadrature, there holds a
sparse tensor product-like error estimate. If \(\varepsilon_\ell\to 0\)
is a monotonically decreasing sequence with \(\varepsilon_\ell\cdot\varepsilon_{L-\ell}
=\varepsilon_L\) for every \(L\in\mathbb{N}\) and
\[
\|\mathcal{Q}_{L-\ell}\rho-\E[\rho]\|\leq c_1\varepsilon_{L-\ell}\quad\text{and}\quad
\|\rho^{(\ell)}-\rho\|\leq c_2\varepsilon_\ell
\]
for some suitable norms and constants \(c_1,c_2>0\), then
\[
\|\mathcal{Q}_{L}^{\text{ML}}[\rho]-\E[\rho]\|\leq C L\varepsilon_L
\]
for a constant \(C>0\). We refer to \cite{HPS12} for the details.

For the calculation of the variance, we employ formula \eqref{eq:cor} and
obtain
\begin{align*}
  \Cor[u_{\mathrm{s}}]({\bs x},{\bs x}') &\approx
	\int_T\int_T
  	\bigg\{\Phi({\bs x}, {\bs z})\overline{\Phi({\bs x}', {\bs z}')}
	\mathcal{Q}_{L}^{\text{ML}}\bigg[\frac{\partial u_{\mathrm{s}}}{\partial{\bs n}}\otimes \frac{\partial u_{\mathrm{s}}}{\partial{\bs n}}\bigg]({\bs z},{\bs z}')\\
  &\phantom{=\int_T\int_T} + \Phi({\bs x}, {\bs z})\overline{\frac{\partial\Phi({\bs x}', {\bs z}')}{\partial {\bs n}_{{\bs z}'}}}
  	\mathcal{Q}_{L}^{\text{ML}}\bigg[\frac{\partial u_{\mathrm{s}}}{\partial{\bs n}}\otimes u_{\mathrm{s}}\bigg]({\bs z},{\bs z}')\\
  &\phantom{=\int_T\int_T} + \frac{\partial\Phi({\bs x}, {\bs z})}{\partial {\bs n}_{\bs z}}\overline{\Phi({\bs x}', {\bs z}')}
  	\mathcal{Q}_{L}^{\text{ML}}\bigg[u_{\mathrm{s}}\otimes\frac{\partial u_{\mathrm{s}}}{\partial{\bs n}}\bigg]({\bs z},{\bs z}')\\
  &\phantom{=\int_T\int_T} + \frac{\partial\Phi({\bs x}, {\bs z})}{\partial {\bs n}_{\bs z}}
	\overline{\frac{\partial\Phi({\bs x}', {\bs z}')}{\partial {\bs n}_{{\bs z}'}}}
	\mathcal{Q}_{L}^{\text{ML}}[u_{\mathrm{s}}\otimes u_{\mathrm{s}}]({\bs z},{\bs z}')\bigg\}\d\sigma_{{\bs z}'}\d\sigma_{\bs z}
\end{align*}
with 
\[
\mathcal{Q}_{L}^{\text{ML}}[\rho\otimes\mu]({\bs z},{\bs z}')\isdef\sum_{\ell=0}^L
\mathcal{Q}_{L-\ell}\big((\rho\otimes\mu)^{(\ell)}({\bs z},{\bs z}',\cdot)
-(\rho\otimes\mu)^{(\ell-1)}({\bs z},{\bs z}',\cdot)\big),
\]
where \((\rho\otimes\mu)^{(\ell)}\isdef\rho^{(\ell)}\otimes\mu^{(\ell)}\).
In principle, it would also be possible to
opt for the multi-index quadrature, which has been proposed
in \cite{BSZ11} for the computation of higher order moments. In this case,
one ends up with
\[
(\rho\otimes\mu)^{(\ell)}({\bs z},{\bs z}',{\bs y})
\isdef\sum_{j=0}^\ell\rho^{(\ell-j)}({\bs z},{\bs y})\mu^{(j)}({\bs z}',{\bs y})
=\sum_{j=0}^\ell\rho^{(j)}({\bs z},{\bs y})\mu^{(\ell-j)}({\bs z}',{\bs y}).
\]
Finally, we remark that there holds a similar error estimate as for the expectation and that isogeometric analysis was recently combined with a multi-index quadrature in \cite{BTT2019}.

\section{Bayesian shape inversion}\label{sec:bayes}
Let $\mathcal{A}({\bs y})\colon H^{1/2}\big(S({\bs y})\big)\to H^{1/2}(T)$, ${\bs y}\in\Gamma$,
be the solution operator which maps the incident wave at \({S}({\bs y})\)
to the scattered wave at \(T\).
Fixing the incident wave \(u_{\text{inc}}\),
we denote by 
\[
G\colon\Gamma\to H^{1/2}(T),\quad {\bs y}\mapsto u_{\text{s}}({\bs y})
\]
the uncertainty-to-solution map.

In forward uncertainty quantification, the goal is to compute quantities
of interest \(\operatorname{QoI}(u_{\text{s}})\), with respect to the
prior measure
$\mu_0$, which is induced by the random variables from \eqref{eq:randomVectorFieldBoundary}.
Often, quantities of interest are assumed to be linear functionals.
The goal of Bayesian inverse uncertainty quantification as in \cite{DS13}
is to incorporate noisy measurements
of solutions \(\mathcal{A}({\bs y})u_{\text{inc}}\),
after potentially incomplete observations.
This is modeled by first considering a bounded, linear
observation operator $O\colon H^{1/2}(T)\to\mathbb{C}^N$, which
models point measurements of the scattered wave at the artificial
interface \(T\). Combining the solution operator with the observation
operatore yields the {uncertainty-to-observation mapping}
\begin{equation}\label{eq:UncToObs}
    \mathcal{G}\colon
    \Gamma\to\mathbb{C}^N
    ,\quad
    {\bs y} \mapsto \mathcal{G}({\bs y}) = O\big(\mathcal{A}({\bs y})u_{\text{inc}}\big)
    .
\end{equation}

The measured data $\bs\delta$ is modeled as resulting from
an observation by $O$, perturbed by additive Gaussian noise
according to
\[
\bs\delta = \mathcal{G}({\bs y}^\star)+{\bs\eta},
\]
where ${\bs y}^\star$ is the unknown, exact parameter. We assume that the noise $\bs\eta$ is given by a complex, circular, symmetric Gaussian random vector with symmetric, positive definite covariance matrix $\bs\Sigma\in\mathbb{R}^{N\times N}$, i.e., $\bs\eta\sim\mathcal{CN}(0,\bs\Sigma)$. Note that this is equivalent to $\bs\eta=\bs\eta_{\text{r}}+i\bs\eta_{\text{i}}$ with uncorrelated  $\bs\eta_{\text{r}}$, $\bs\eta_{\text{i}}$ and $\bs\eta_{\text{r}},\bs\eta_{\text{i}}\sim\mathcal{N}(0,\Sigma/2)$, and respects the physical time-harmonic model of the scattering problem, see \cite{TV05}.

Within this article, we aim at predicting the shape of the
random scatterer based on observations of \(u_{\text{s}}\) at \(T\).
Concretely, we wish to compute expectation and variance of the
deformation field.
To that end, we define the Gaussian potential, also referred to
as the least-squares or data misfit functional, by
$\Phi_{\bs\Sigma}\colon \Gamma\times\mathbb{C}^N\to\mathbb{R}$,
\begin{equation}\label{eq:pot}
    \Phi_{\bs\Sigma}({\bs y},\bs\delta)\isdef
    \frac12 \|\bs\delta-\mathcal{G}({\bs y})\|^2_{\bs\Sigma}
    \isdef \frac12 \big(\bs\delta-\mathcal{G}({\bs y})\big)^\intercal {\bs\Sigma}^{-1} \big(\bs\delta-\mathcal{G}({\bs y})\big)
    .
\end{equation}

Given the prior measure $\mu_0$,
Bayes' formula yields an expression for the posterior measure
$\mu^\delta$ on $\Gamma$, given the data $\bs\delta$ with the 
Radon-Nikodym derivative is given by
\[
    \frac{d\mu^{\bs\delta}}{d\mu_0}({\bs y}) 
    =
    \frac{e^{ -\Phi_{\bs\Sigma}({\bs y},{\bs\delta})}}{Z}  
\]
with
\[Z\isdef \int_\Gamma e^{
-\Phi_{\bs\Sigma}({\bs y},\bs\delta)} \,\mu_0(\d{\bs y}) > 0,
\] 
see \cite{DS13}.

Now, the expected shape of the random scatterer is given by
\[
\E^{\mu_{\bs\delta}}[\bs\chi](S_\refd)
\isdef\int_\Gamma\bs\chi(S_\refd,{\bs y})\frac{e^{ -\Phi_{\bs\Sigma}({\bs y},{\bs\delta})}}{Z}
\d\mu_0({\bs y})
\]
and its variance by
\[
\V^{\mu_{\bs\delta}}[\bs\chi](S_\refd)
\isdef\int_\Gamma\bs
\chi(S_\refd,{\bs y})\bs\chi(S_\refd,{\bs y})^\intercal\frac{e^{ -\Phi_{\bs\Sigma}({\bs y},{\bs\delta})}}{Z}
\d\mu_0({\bs y})-\E^{\mu_{\bs\delta}}[\bs\chi](S_\refd)\E^{\mu_{\bs\delta}}[\bs\chi](S_\refd)^\intercal.
\]

In order to approximate these integrals numerically, we shall
employ the multilevel ratio estimator, which splits the computation
of the actual integral and the normalization constant and approximates
each by a telescoping sum, see \cite{DGLS17} and the references therein. For the normalization constant,
we consider
\[
\mathcal{Q}_{L}^{\text{ML}}[\rho]\isdef\sum_{\ell=0}^L
\mathcal{Q}_{L-\ell}(\rho_\ell-\rho_{\ell-1})
\]
with
\[
\rho_\ell\isdef e^{
-\Phi_{{\bs\Sigma},\ell}({\bs y},\bs\delta)},\quad
\Phi_{{\bs\Sigma},\ell}({\bs y},{\bs\delta})\isdef\big\|
\bs\delta-O\big(\mathcal{A}_\ell({\bs y})u_{\text{inc}}\big)\big\|_{\bs\Sigma},
\quad Z_{-1}\isdef 0,
\]
i.e.\ we consider a multilevel hierarchy based on approximations of
the scattered wave on different levels of refinement.
Now, we may compute, for example, the expected deformation field
according to
\[
\mathcal{Q}_{L}^{\text{ML},\mu_{\bs\delta}}[\bs\chi]\isdef\bigg(\sum_{\ell=0}^L
\mathcal{Q}_{L-\ell}\big({\bs\chi}\cdot(\rho_\ell-\rho_{\ell-1})\big)
\bigg)\bigg/\mathcal{Q}_{L}^{\text{ML}}[\rho].
\]

\section{Numerical examples}\label{sec:numex}
\subsection{Geometries, discretization, and multilevel quadrature}
\begin{figure}
\centering
\includegraphics[width=0.46\textwidth,clip=true,trim=90 140 350 350]{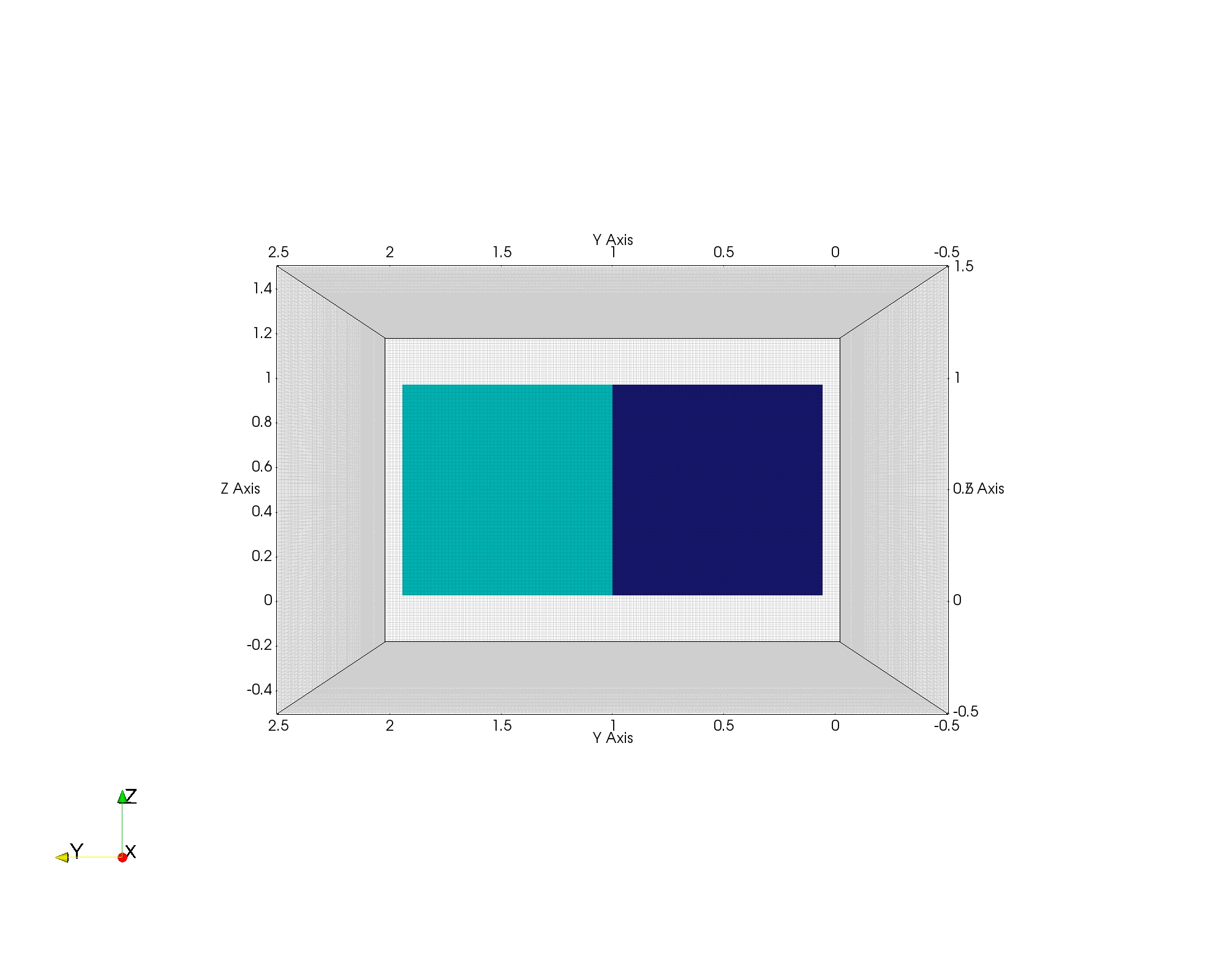}\quad
\includegraphics[width=0.46\textwidth,clip=true,trim=90 140 220 380]{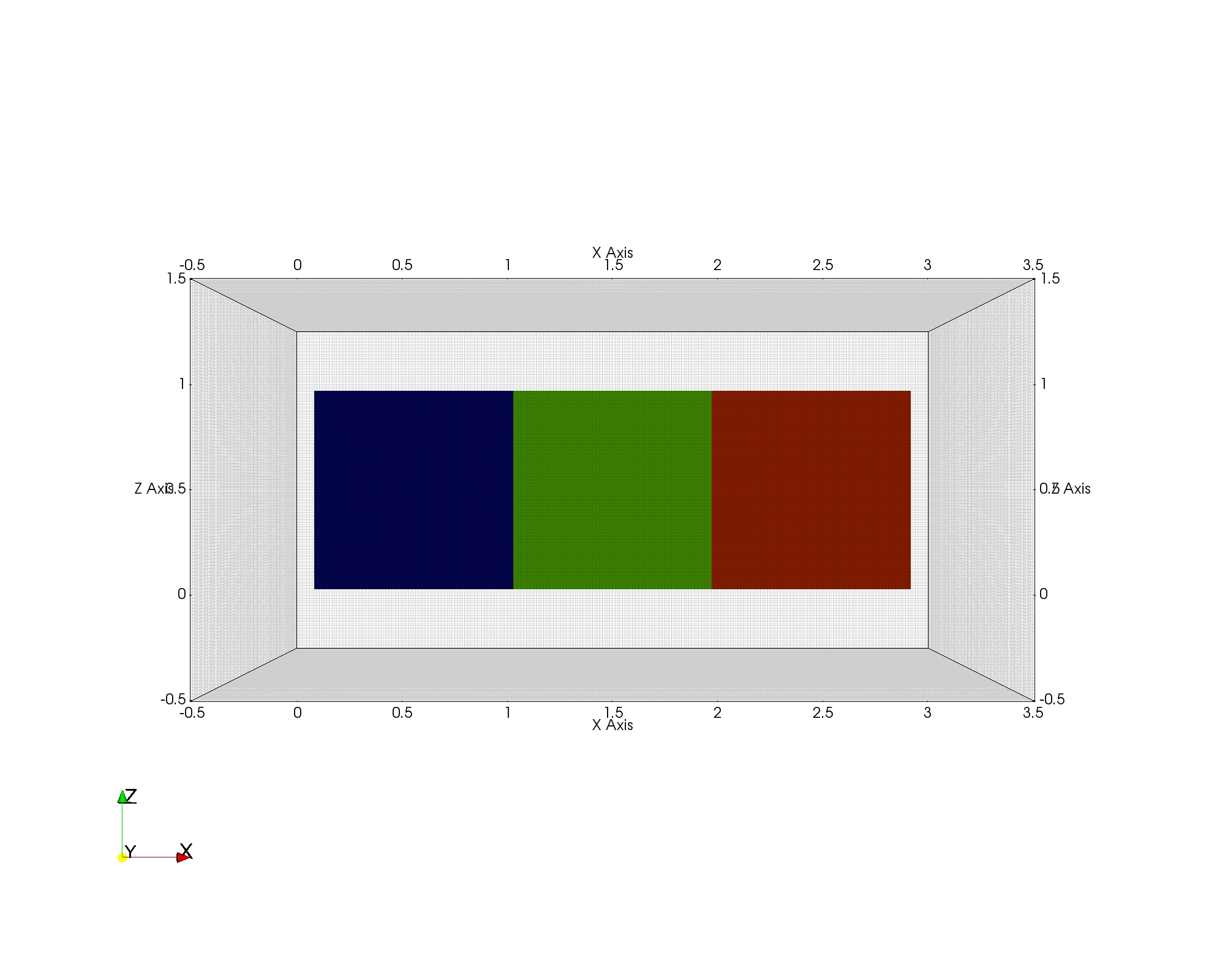}
\includegraphics[width=0.46\textwidth,clip=true,trim=90 140 250 250]{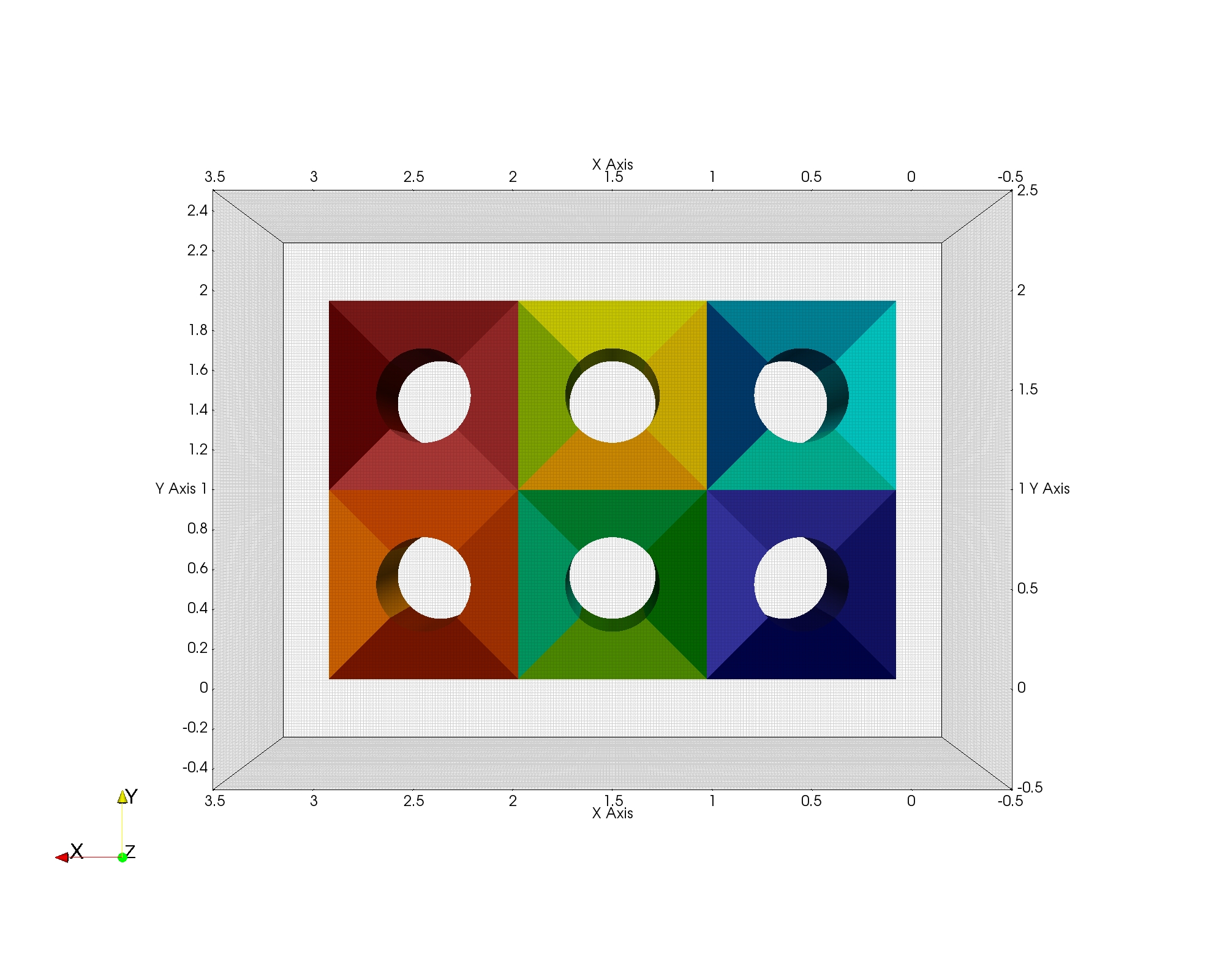}
\caption{The $[0,3]\times[0,2]\times[0,1]$ cuboid with drilled holes and the $[-1.5,3.5]\times[-0.5,2.5]\times[-0.5,1.5]$ artificial interface.}
\label{fig:cuboids}
\end{figure}

We consider a scatterer $D_\refd$ given by a cuboid $[0,3]\times[0,2]\times[0,1]$ with six drilled holes,
with an artificial interface $T$ given by the cuboid
$[-1.5,3.5]\times[-0.5,2.5]\times[-0.5,1.5]$. A visualization of the situation
may be found in Figure~\ref{fig:cuboids}.
The surface of the scatterer is represented by 82 patches as illustrated, see again Figure~\ref{fig:cuboids}, and the artificial interface by 52 patches. The wavenumber is chosen as $\kappa=1$.

We discretize the random field with globally continuous B-splines of polynomial
degree $p=2$ in each spatial variable and uniform three spatial refinements,
leading to a dense covariance matrix
$\underline{\bfC}\in\mathbb{R}^{19\,896\times 19\,896}$.
For an efficient computation of the Karhunen-Lo\`eve expansion, we proceed as
outlined in Remark~\ref{rem:efficientIGAKL}. The artificial interface is discretized
with tensor-product polynomials of degree $p=6$ on each patch. The Cauchy data on the artificial interface can then be obtained from the values on $52\cdot 7^2=2\,548$ point evaluations on the interface by solving $52$ local interpolation problems of size $7^2$.

For the application of the multilevel quadrature, we perform the acoustic scattering computations with patchwise continuous B-splines of degree $p=2$ and the refinement levels $\ell=0,1,2,3$, leading to 738, 1\,312, 2\,952, and 8\,200 degrees of freedom. The implementation of the spatial discretizations is based on the \verb|C++|-library \verb+bembel+ \cite{bembel, DHK+20}, which is easily adapted to our needs and provides fast compression schemes for the scattering computations.

The multilevel quadrature is either based on a quasi-Monte Carlo quadrature using the Halton sequence, see \cite{Caf98}, or on the anisotropic sparse grid quadrature using Gau{\ss}-Legendre points as described in \cite{HHPS18}. The latter is available as open source software package \verb+SPQR+\footnote{https://github.com/muchip/SPQR}. Due to the high asymptotic convergence rate of $h^{2p+2}\sim 2^{-6}$ of the higher-order method for the scattering computations, the number of samples for the multilevel quadrature has to be adapted for each level as shown in Table~\ref{tab:numberofsamples}, where `QMC' stands for the quasi-Monte Carlo quadrature and `SG' for the sparse grid quadrature.

\begin{table}[htb]
    \centering
\begin{tabular}{|c|c|c|c|c|}
    \hline
     & \(\ell=0\) & \(\ell=1\) & \(\ell=2\) &\(\ell=3\)\\\hline
    QMC & 2\,097\,152& 32\,768& 512& 256 \\\hline
    SG  & 2\,328\,341 & 58\,251 & 957 & 351 \\\hline
    \end{tabular}
    \medskip
    \caption{Number of samples on the different levels for SG and QMC.}
    \label{tab:numberofsamples}
\end{table}

Due to the large number of samples and computational costs of solving three dimensional scattering problems, we employ a hybrid parallelization with MPI and OpenMP to accelerate the sampling process. The computations have been carried out with up to 24 MPI processes consisting of up to 8 OpenMP threads each, resulting in a total of up to 192 cores. The OpenMP threads have been dedicated to the boundary element solver, while the MPI processes were used to parallelize the sampling of the random parameter.

\subsection{Forward problem}
\begin{figure}
\begin{center}
\begin{tikzpicture}
\draw(0,0) node{\includegraphics[width=0.38\textwidth,clip=true,trim=300 80 350 250]{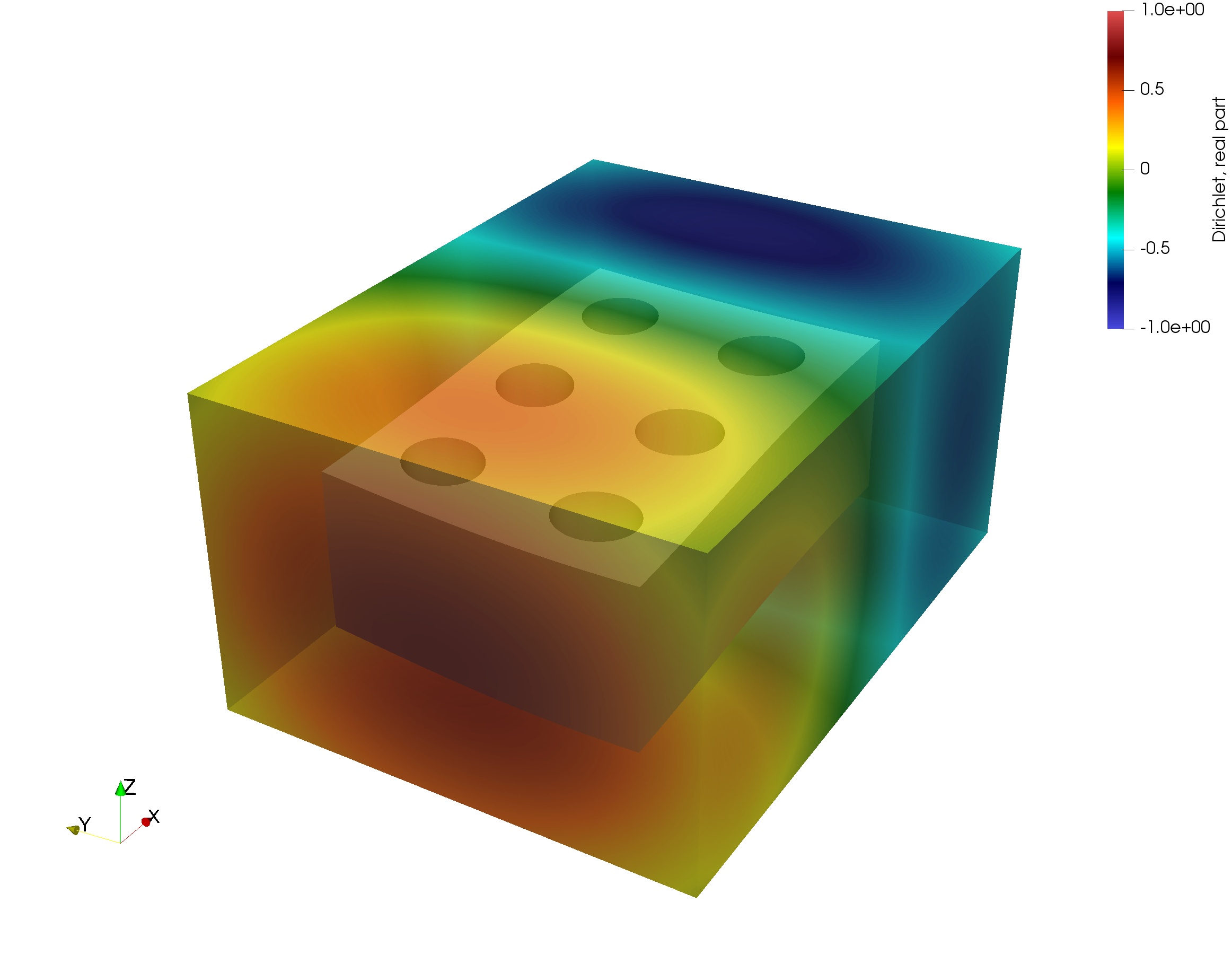}};
\draw(5.8,0) node{\includegraphics[width=0.38\textwidth,clip=true,trim=300 80 350 250]{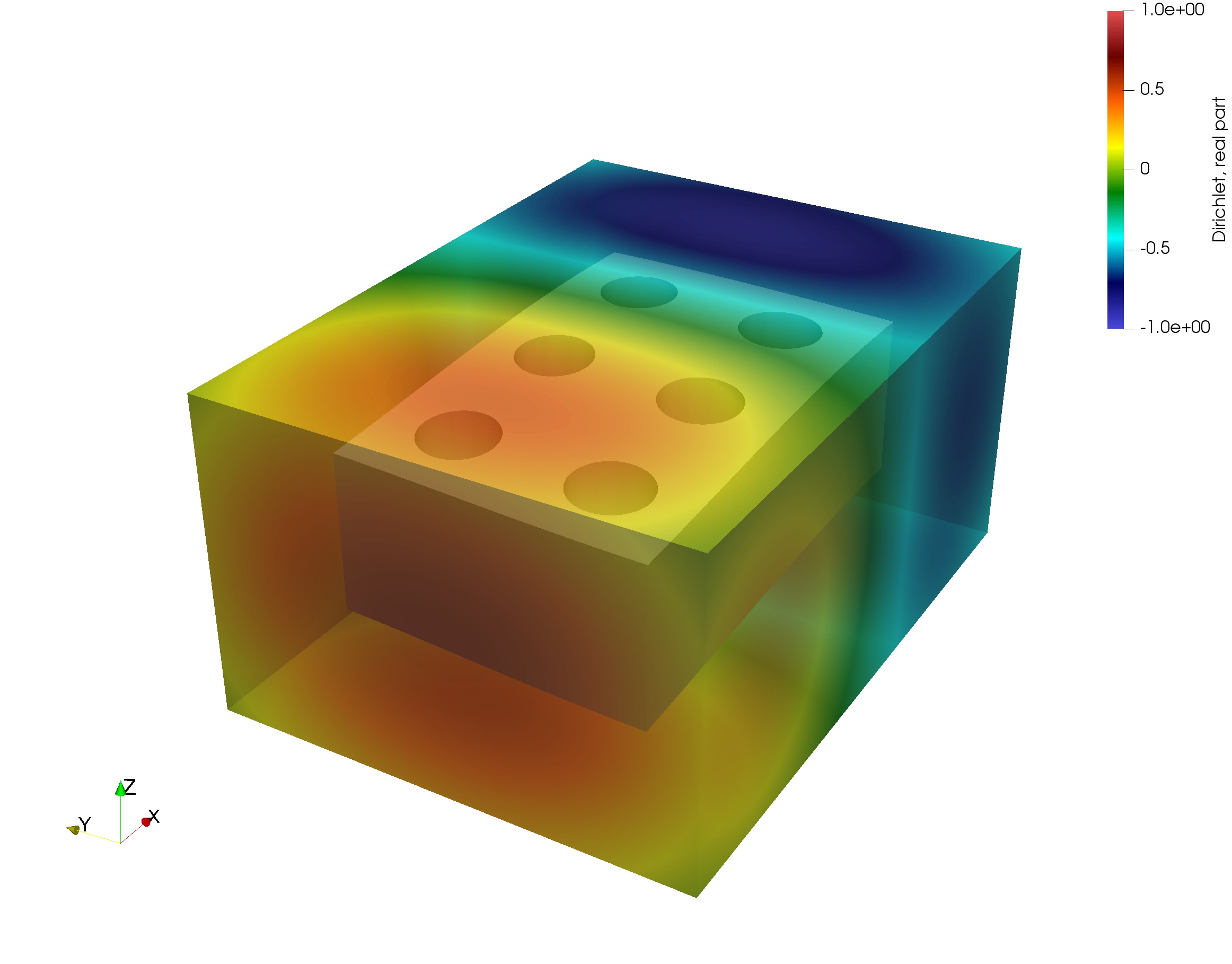}};
\draw(0,5.8) node{\includegraphics[width=0.38\textwidth,clip=true,trim=300 80 350 250]{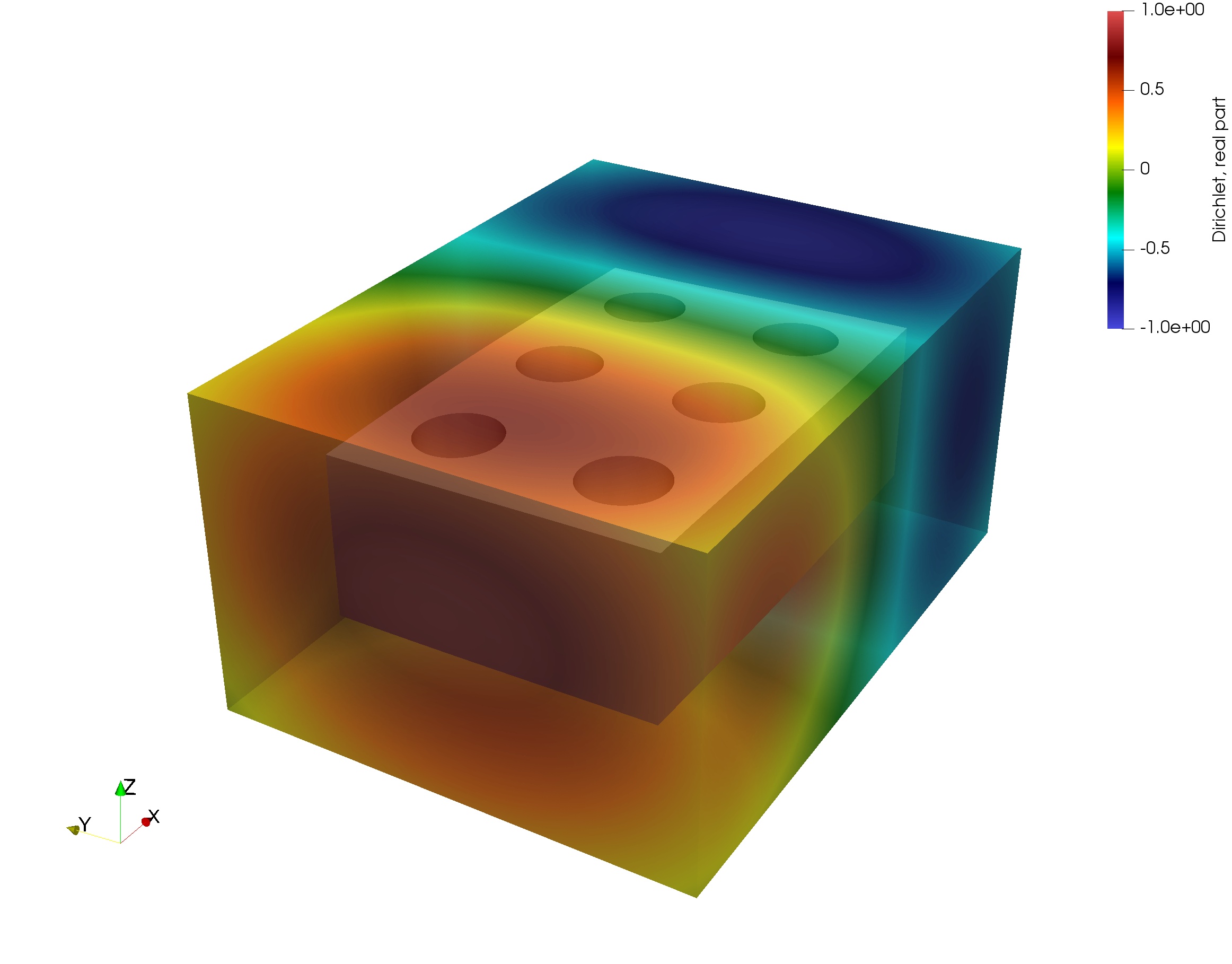}};
\draw(5.8,5.8) node{\includegraphics[width=0.38\textwidth,clip=true,trim=300 80 350 250]{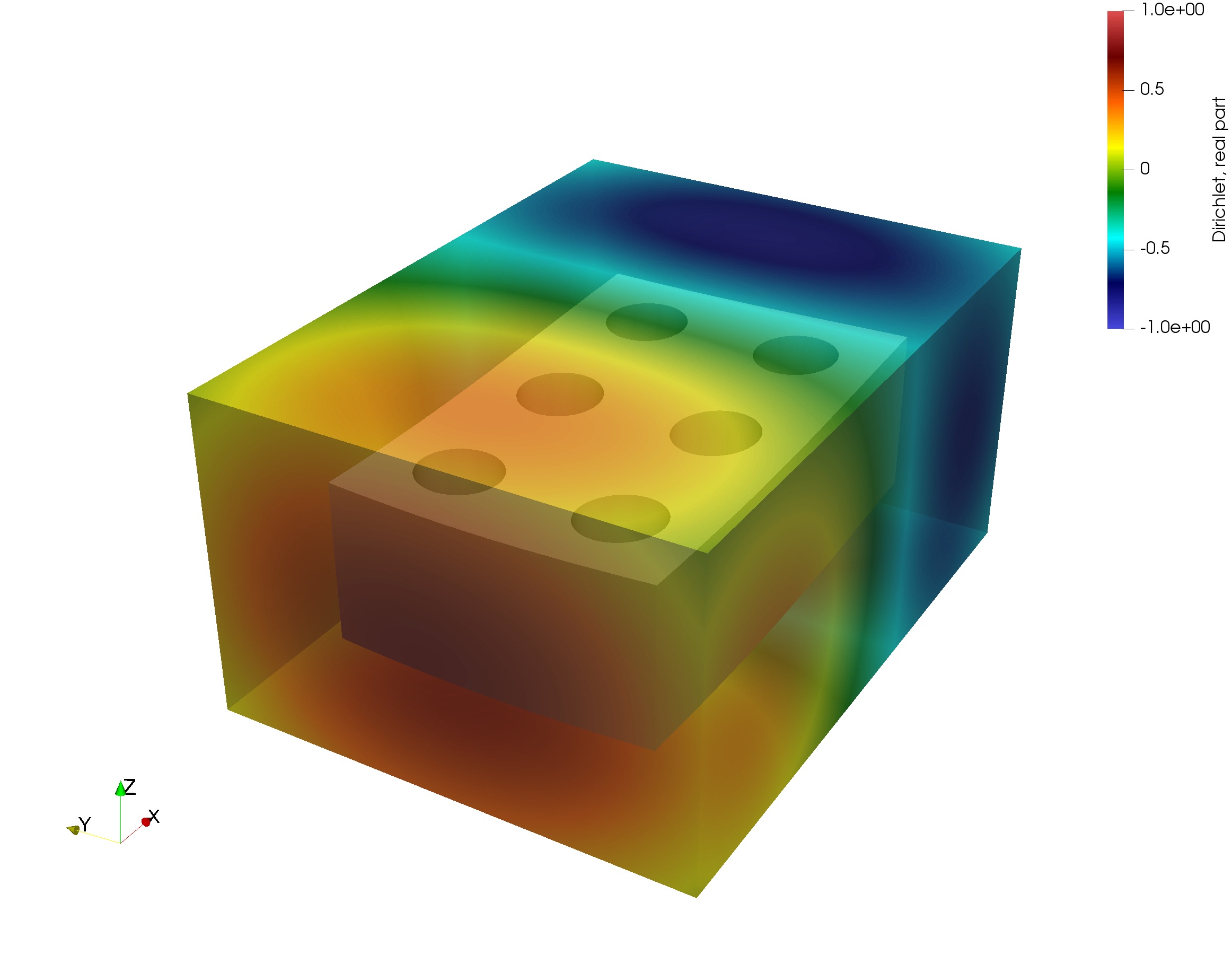}};
\draw(10.2,2.9) node{\includegraphics[width=0.14\textwidth,clip=true,trim=2100 1180 0 0]{samples0002}};
\end{tikzpicture}
\caption{Domain perturbations drawn from the random field and scattered wave on the artificial interface.}
\label{fig:perturbed}
\end{center}
\end{figure}

\begin{figure}
    \centering
\begin{tikzpicture}
\pgfplotstableread{%
  1     0.891113
  2     0.891113
  3     0.891113
  4      0.52293
  5      0.52293
  6     0.522929
  7     0.400666
  8     0.400666
  9     0.400666
 10     0.235048
 11     0.235048
 12     0.235048
 13     0.229026
 14     0.229026
 15     0.229026
 16     0.217782
 17     0.217782
 18     0.217782
 19     0.132593
 20     0.132593
 21     0.132593
 22     0.110386
 23     0.110386
 24     0.110386
 25    0.0998649
 26    0.0998649
 27    0.0998649
 28    0.0942553
 29    0.0942553
 30    0.0942553
 31    0.0707209
 32    0.0707209
 33    0.0707209
 34    0.0618501
 35    0.0618501
 36    0.0618501
 37    0.0556918
 38    0.0556918
 39    0.0556918
 40       0.0546
 41       0.0546
 42       0.0546
 43    0.0308939
 44    0.0308939
 45    0.0308936
 46    0.0302369
 47    0.0302369
 48    0.0302369
 49    0.0280594
 50    0.0280594
 51    0.0280594
 52    0.0260125
 53    0.0260125
 54    0.0260125
 55    0.0249283
 56    0.0249283
 57    0.0249281
 58    0.0227158
 59    0.0227158
 60    0.0227158
 61    0.0197408
 62    0.0197408
 63    0.0197405
 64    0.0178012
 65    0.0178012
 66    0.0178011
 67    0.0177589
 68    0.0177589
 69    0.0177588
 70     0.015502
 71     0.015502
 72     0.015502
 73    0.0141151
 74    0.0141151
 75    0.0141149
 76    0.0135702
 77    0.0135702
 78    0.0135702
 79    0.0086055
 80    0.0086055
 81   0.00860491
 82   0.00836435
 83   0.00836435
 84   0.00836399
 85   0.00764791
 86   0.00764791
 87   0.00764789
 88   0.00737144
 89   0.00737144
 90   0.00737132
 91   0.00734445
 92   0.00734445
 93   0.00734444
 94   0.00632728
 95   0.00632728
 96   0.00632728
 97   0.00595073
 98   0.00595073
 99    0.0059507
100   0.00576797
101   0.00576797
102   0.00576778
103   0.00477007
104   0.00477007
105   0.00477007
106   0.00465568
107   0.00465568
108   0.00465099
109   0.00453527
110   0.00453527
111   0.00452787
112    0.0033021
113    0.0033021
114    0.0033021
115   0.00319293
116   0.00319293
117   0.00319282
118   0.00312031
119   0.00312031
120   0.00312017
121   0.00296512
122   0.00296512
123   0.00296498
124    0.0025919
125    0.0025919
126    0.0025809
127   0.00233408
128   0.00233408
129   0.00233393
130   0.00224389
131   0.00224389
132   0.00224255
133    0.0020813
134    0.0020813
135   0.00208072
136   0.00203516
137   0.00203516
138   0.00203509
139   0.00200514
140   0.00200514
141   0.00200514
142   0.00183949
143   0.00183949
144   0.00183921
145   0.00177817
146   0.00177817
147   0.00177791
148   0.00176621
149   0.00176621
150   0.00176621
151   0.00137309
152   0.00137309
153   0.00137308
154   0.00132827
155   0.00132827
156    0.0013282
157   0.00112145
158   0.00112145
159   0.00112019
160   0.00108044
161   0.00108044
162   0.00108043
163   0.00105601
164   0.00105601
165   0.00105433
                       }
   \mytable;
\begin{semilogyaxis}[xmin=0,xmax=166,ymin=0.9e-3,ymax=1.1,width=0.8\textwidth,height=0.5\textwidth,grid=both]
\addplot[mark=*,mark size=1pt, only marks]table{\mytable};
\end{semilogyaxis}
\end{tikzpicture}
    \caption{Numerical approximation of the singular values of the covariance operator under consideration.}
    \label{fig:singularvalues}
\end{figure}
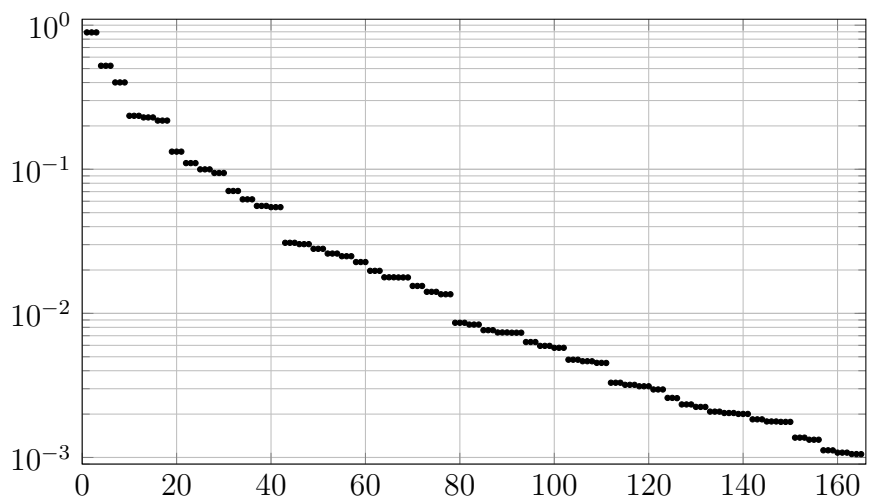

We consider a centered Gaussian random field for the domain perturbations with covariance function
\[
\Cov[\boldsymbol{\chi}_S](\bfx,{\bs y})
=
\frac{1}{20}
\begin{bmatrix}
e^{-\frac{\|\bfx-{\bs y}\|_2^2}{4}} & 0 & 0 \\
0 & e^{-\frac{\|\bfx-{\bs y}\|_2^2}{4}} & 0 \\
0 & 0 & e^{-\frac{\|\bfx-{\bs y}\|_2^2}{4}}
\end{bmatrix}.
\]
Four different realizations of the deformed scatterer and corresponding scattered waves at the artificial interface are illustrated in Figure~\ref{fig:perturbed}. The singular values of the corresponding deformation field are illustrated in Figure~\ref{fig:singularvalues}. The parameter dimension
is 165.

For demonstrating the validity of the dimension reduction via the artificial interface, we also define 100 evaluation points outside of the artificial interface, which are equally distributed on a sphere centered around the origin with radius $5$. Note that the origin is one of the corners of the reference geometry.

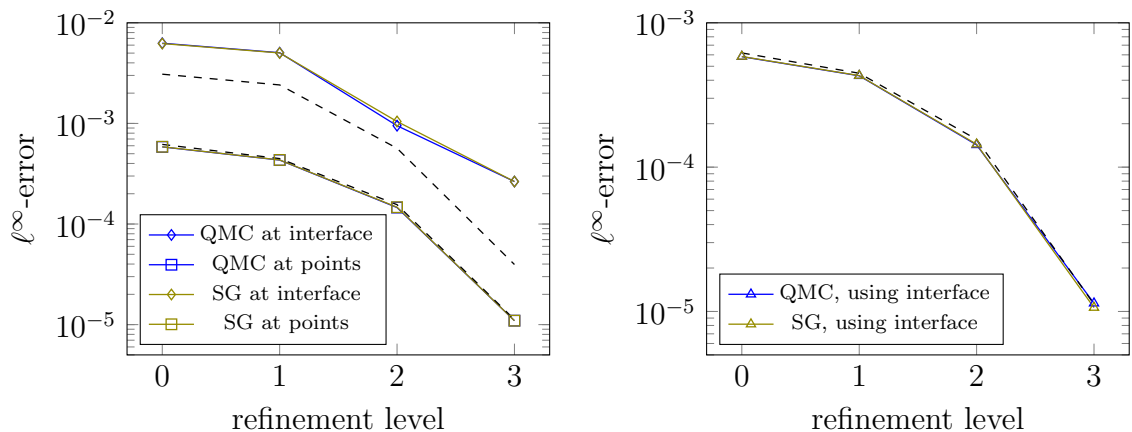
\begin{figure}[htb]
\begin{center}
\begin{tikzpicture}
\begin{semilogyaxis}[height=0.4\textwidth,width=0.48\textwidth, ymax = 1e-2, ymin = 5e-6,
    legend style={legend pos=south west,font=\tiny}, xlabel=refinement level, ylabel ={$\ell^{\infty}$-error}]
\addplot [line width=0.5pt, color=blue,mark=diamond] table[x index=0,y index=1]{./data/MLQMC_mean_interface.txt};\addlegendentry{QMC at interface};
\addplot [line width=0.5pt, color=blue,mark=square] table[x index=0,y index=1]{./data/MLQMC_mean_points.txt};\addlegendentry{QMC at points};
\addplot [line width=0.5pt, color=olive,mark=diamond] table[x index=0,y index=1]{./data/MLSG_mean_interface.txt};\addlegendentry{SG at interface};
\addplot [line width=0.5pt, color=olive,mark=square] table[x index=0,y index=1]{./data/MLSG_mean_points.txt};\addlegendentry{SG at points};
\addplot[line width=0.5pt, color=black,mark=none, dashed] coordinates {
(0,0.00308649)
(1,0.00241284)
(2,0.000565987)
(3,3.95972e-05)
};
\addplot[line width=0.5pt, color=black,mark=none, dashed] coordinates {
(0,0.000618693)
(1,0.000447365)
(2,0.000155738)
(3,1.12308e-05)
};
\end{semilogyaxis}
\end{tikzpicture}
\hfill
\begin{tikzpicture}
\begin{semilogyaxis}[height=0.4\textwidth,width=0.48\textwidth, ymax = 1e-3, ymin = 5e-6,
    legend style={legend pos=south west,font=\tiny}, xlabel=refinement level, ylabel ={$\ell^{\infty}$-error}]
\addplot [line width=0.5pt, color=blue,mark=triangle] table[x index=0,y index=2]{./data/MLQMC_mean_points.txt};\addlegendentry{QMC, using interface};
\addplot [line width=0.5pt, color=olive,mark=triangle] table[x index=0,y index=2]{./data/MLSG_mean_points.txt};\addlegendentry{SG, using interface};
\addplot[line width=0.5pt, color=black,mark=none, dashed] coordinates {
(0,0.000618693)
(1,0.000447365)
(2,0.000155738)
(3,1.12308e-05)
};
\end{semilogyaxis}
\end{tikzpicture}
\caption{\label{fig:QuadProb4}{\em Left:} Convergence of the multilevel quadrature against the MLQMC solution for the mean over the artificial interface and on points on the sphere. {\em Right:} Convergence in the points on the sphere when they are evaluated from the mean of the Cauchy data at the artificial interface.}
\end{center}
\end{figure}

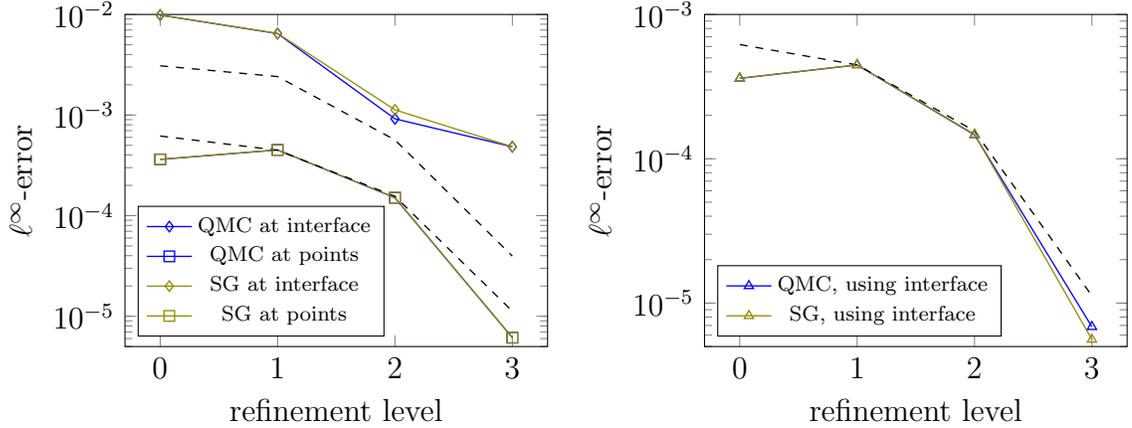
\begin{figure}[htb]
\begin{center}
\begin{tikzpicture}
\begin{semilogyaxis}[height=0.4\textwidth,width=0.48\textwidth, ymax = 1e-2, ymin =5e-6,
    legend style={legend pos=south west,font=\tiny}, xlabel=refinement level, ylabel ={$\ell^{\infty}$-error}]
\addplot [line width=0.5pt, color=blue,mark=diamond] table[x index=0,y index=1]{./data/MLQMC_cor_interface.txt};\addlegendentry{QMC at interface};
\addplot [line width=0.5pt, color=blue,mark=square] table[x index=0,y index=1]{./data/MLQMC_cor_points.txt};\addlegendentry{QMC at points};
\addplot [line width=0.5pt, color=olive,mark=diamond] table[x index=0,y index=1]{./data/MLSG_cor_interface.txt};\addlegendentry{SG at interface};
\addplot [line width=0.5pt, color=olive,mark=square] table[x index=0,y index=1]{./data/MLSG_cor_points.txt};\addlegendentry{SG at points};
\addplot[line width=0.5pt, color=black,mark=none, dashed] coordinates {
(0,0.00308649)
(1,0.00241284)
(2,0.000565987)
(3,3.95972e-05)
};
\addplot[line width=0.5pt, color=black,mark=none, dashed] coordinates {
(0,0.000618693)
(1,0.000447365)
(2,0.000155738)
(3,1.12308e-05)
};
\end{semilogyaxis}
\end{tikzpicture}\hfill
\begin{tikzpicture}
\begin{semilogyaxis}[height=0.4\textwidth,width=0.48\textwidth, ymax = 1e-3, ymin = 5e-6,
    legend style={legend pos=south west,font=\tiny}, xlabel=refinement level, ylabel ={$\ell^{\infty}$-error}]
\addplot [line width=0.5pt, color=blue,mark=triangle] table[x index=0,y index=2]{./data/MLQMC_cor_points.txt};\addlegendentry{QMC, using interface};
\addplot [line width=0.5pt, color=olive,mark=triangle] table[x index=0,y index=2]{./data/MLSG_cor_points.txt};\addlegendentry{SG, using interface};
\addplot[line width=0.5pt, color=black,mark=none, dashed] coordinates {
(0,0.000618693)
(1,0.000447365)
(2,0.000155738)
(3,1.12308e-05)
};
\end{semilogyaxis}
\end{tikzpicture}
\caption{\label{fig:QuadProb5}{\em Left:} Convergence of the multilevel quadrature against the MLQMC solution for the correlation in the points at the interface and in the points in free space. {\em Right:} Convergence in the points in free space when computing the correlation in free space from the correlation at the artificial interface.}
\end{center}
\end{figure}

In order to measure approximation errors, we compare the solutions obtained by the multilevel sparse grid quadrature to that of the multilevel quasi-Monte Carlo
quadrature on the finest level \(L=3\) and vice versa. The left-hand side of Figure~\ref{fig:QuadProb4} shows the convergence error of the expectation for the Cauchy data at the interface and for the 100 points on the sphere. The dashed curves indicate the convergence of the spatial approximation on the reference domain. The right-hand side of Figure~\ref{fig:QuadProb4} illustrates the convergence of the potential evaluation on the sphere when using the mean of the Cauchy data on the artificial interface, i.e., when using \eqref{eq:MLQmean}.
Figure~\ref{fig:QuadProb5} shows these quantities for the correlation.

\subsection{Shape inversion}
\begin{figure}
\centering
\includegraphics[width=0.46\textwidth,clip=true,trim=90 140 350 350]{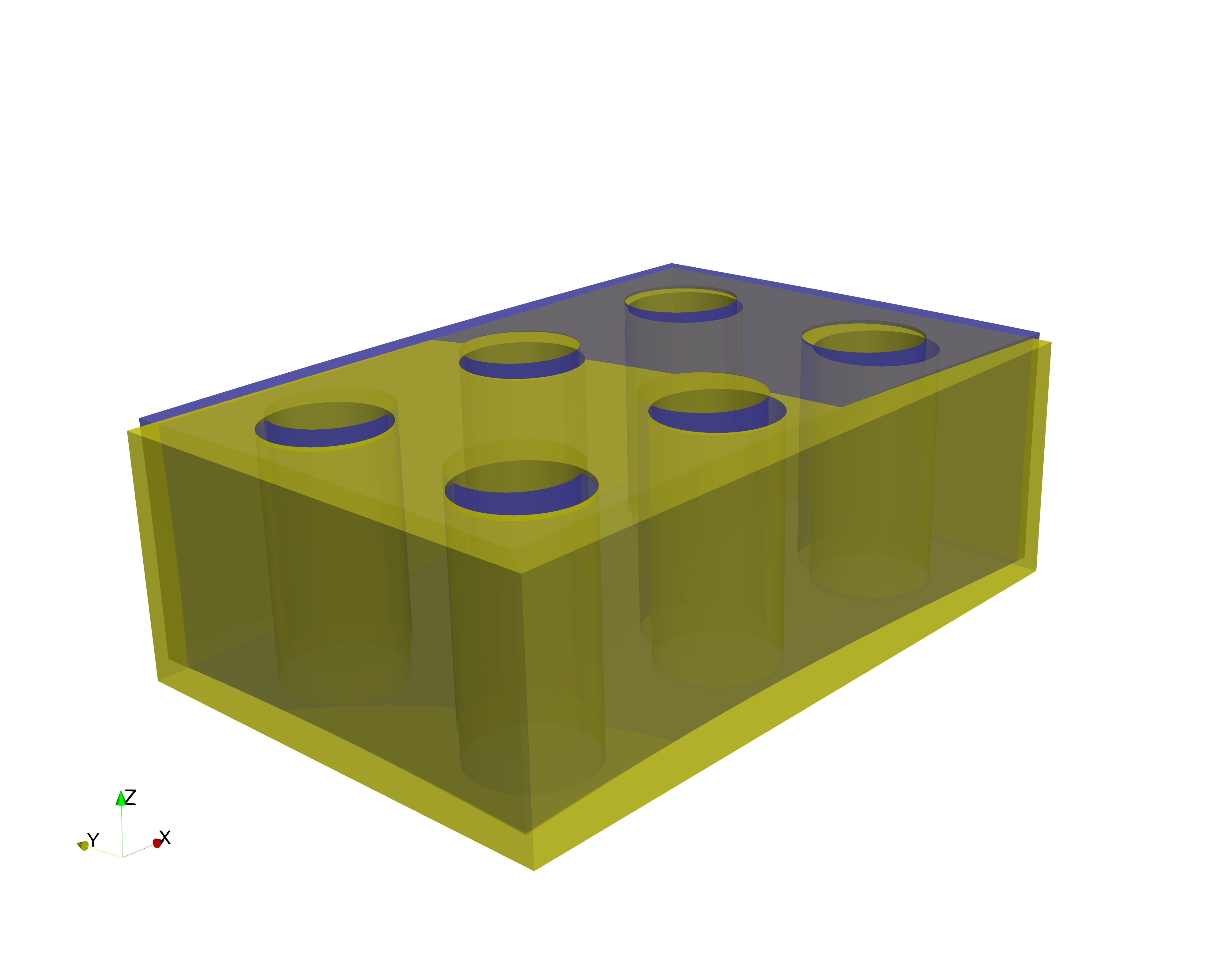}
\includegraphics[width=0.46\textwidth,clip=true,trim=90 140 350 350]{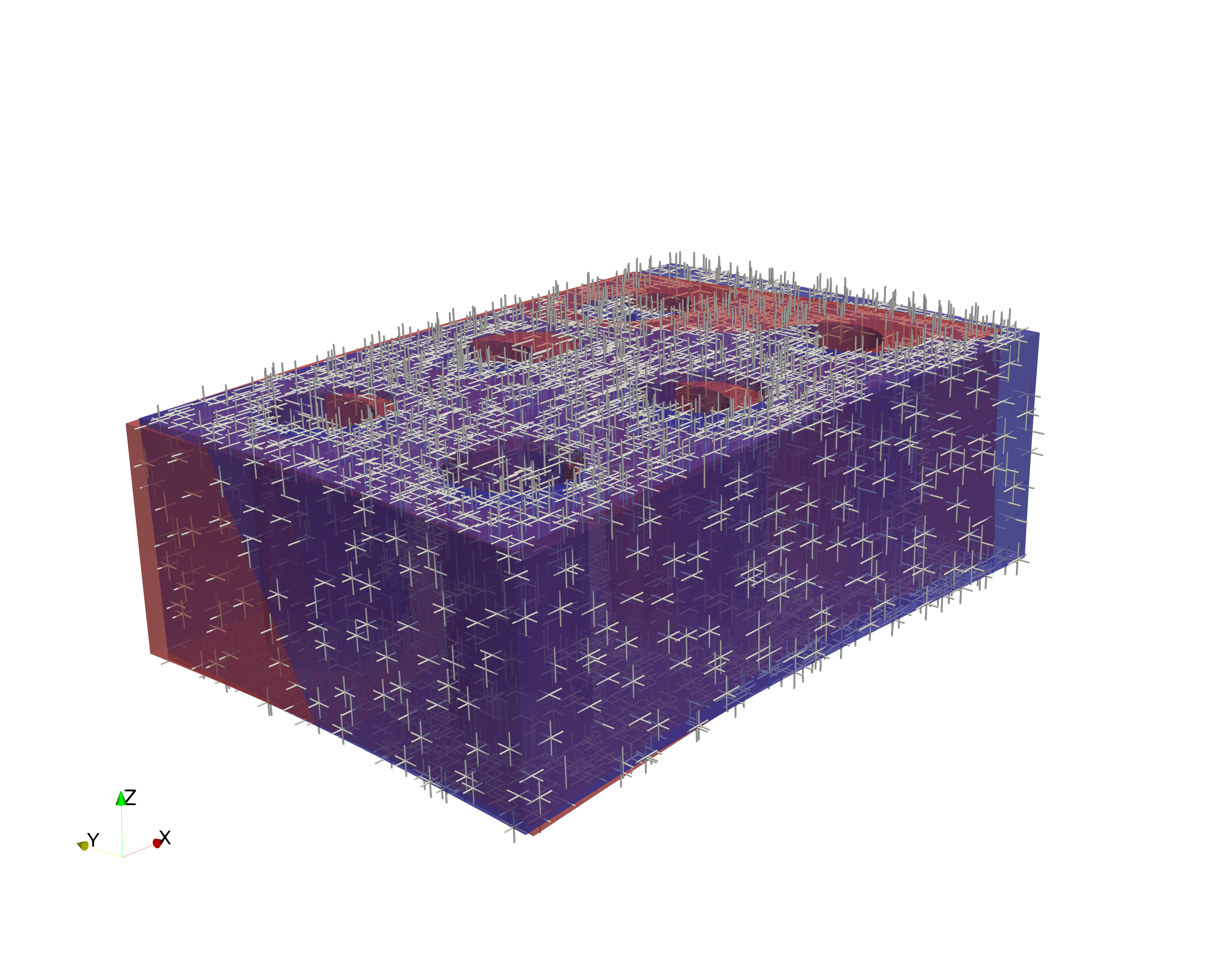}
\caption{\emph{Left:} Prior mean (yellow) and posterior mean (blue) of the inverse problem. \emph{Right:} Reference solution (red) and posterior mean (blue) of the inverse problem with $2\sigma$ confidence intervals (crosses) in each coordinate direction.}
\label{fig:bayes}
\end{figure}

For illustrating the Bayesian shape inversion, we draw a random domain perturbation given by ${\bs y}^\star\in\Gamma$ from the model presented in the previous subsection and consider it to be our reference solution. The measurement operator $O$ defining $\mathcal{G}$ is given by point evaluations of the scattered wave in the midpoints of the 52 patches at the artificial interface. The noise level is set to $\bs\Sigma=\sigma^2{\bs I}$, where $\sigma=0.1\cdot\max |\mathcal{G}({\bs y}^\star)|$. The unperturbed domain is considered as a prior. Figure \ref{fig:bayes} illustrates the reference solution, the prior and posterior mean and the posterior's $2\sigma$ confidence region in each coordinate direction obtained by a multilevel quasi-Monte Carlo quadrature on the finest level \(L=3\). The posterior mean has clearly moved away from the prior and is located within the $2\sigma$ region of the true scatterer.

\section{Conclusions and Future work}

We have introduced a fast IGA implementation for solving time-harmonic acoustic wave scattering for shape uncertainty quantification, employing boundary integral formulations, multi-level quadrature and state-of-the art acceleration techniques. This allows for the analysis of large shape deformations for both forward and inverse problems, including shape optimization. Future work involves the extension to Maxwell scattering.

\bibliographystyle{plain}
\bibliography{literature}
\end{document}